\documentclass[11pt]{article}

\title{
Rate-Induced Tipping in Heterogeneous Reaction-Diffusion Systems:\\ An Invariant Manifold Framework and  Geographically Shifting Ecosystems
}

\author{Cris R. Hasan\thanks{Adam Smith Business School, University of Glasgow, Glasgow, UK (cris.hasan@glasgow.ac.uk).} \and Ruaidhrí Mac Cárthaigh\thanks{School of Mathematical Sciences, University College Cork, Cork, Ireland.} \and Sebastian Wieczorek\footnotemark[2]}

\date{}

\usepackage[utf8]{inputenc}
\usepackage{comment}
\usepackage{graphicx}
\usepackage{amssymb}
\usepackage{amsmath}
\usepackage{hyperref}
\usepackage[numbers]{natbib}
\usepackage{cleveref}
\usepackage{amsthm}
\usepackage[font=small,labelfont=bf]{caption}
\usepackage[dvipsnames]{xcolor}
\usepackage{todonotes}
\usepackage[left=2cm]{geometry}
\usepackage[normalem]{ulem}
\definecolor{mygray}{gray}{0.6}
\usepackage{url}

\DeclareMathOperator{\sech}{sech}

\newcommand{\ch}{\color{red}}

\numberwithin{equation}{section}

\hypersetup{
    colorlinks,
    linkcolor={red!70!black},
    citecolor={blue!70!black},
    urlcolor={blue!80!black}
}

\begin{document}
\maketitle
\section*{Abstract}
We propose a framework to study tipping points in reaction-diffusion equations (RDEs) in one spatial dimension, where the reaction term decays in space (asymptotically homogeneous) and varies linearly with time (nonautonomous) due to an external input. 
A compactification of the moving-frame coordinate together with Lin's method to construct heteroclinic orbits along intersections of stable and unstable invariant manifolds allows us to (i) obtain multiple coexisting  pulse and front solutions for the RDE by computing heteroclinic orbits connecting equilibria at negative and positive infinity  in the compactified  moving-frame ordinary differential equation, 
(ii) detect tipping points as dangerous bifurcations of such heteroclinic orbits and, (iii) obtain tipping diagrams by numerical continuation of such bifurcations.
We apply our framework to an illustrative model  of a habitat patch that features an Allee effect in population growth and is geographically shrinking or shifting due to human activity or climate change.
Thus, we identify two classes of tipping points 
to extinction: bifurcation-induced tipping (B-tipping) when the shrinking habitat falls below  some critical length and  rate-induced tipping (R-tipping) when the shifting habitat exceeds some critical speed.
We explore two-parameter R-tipping diagrams to understand how the critical speed depends on the size of the habitat patch and the dispersal rate of the population, uncover parameter regions where the shifting population survives, and relate 
these regions to the invasion speed in an infinite homogeneous habitat. 
Furthermore, we contrast the tipping instabilities  with gradual transitions to extinction found for  logistic population growth  without the Allee effect.
\subsubsection*{Keywords}
Tipping points, reaction-diffusion equations, moving habitats, compactification, numerical continuation, invariant manifolds, Lin's method, regime shifts.

\newpage
\tableofcontents
\section{Introduction}
Tipping points or critical transitions
are often described as large and sudden changes in the state of an open system that arise in response to small and slow changes in the external inputs.
The phenomenon of tipping is ubiquitous in natural and human systems, could be of great environmental impact, and  has thus attracted much interest from the scientific community over the past two decades, especially in climate science~\cite{lenton2011early,lenton2019climate}, as well as in ecology~\cite{
%abbot2011bifurcations,hill2016analysis,lenton2011early,lenton2019climate,scheffer2001catastrophic,
gowda2014transitions,nobre2009tipping,scheffer2009early,selkoe2015principles,van2016you}, where it is referred to as a ``regime shift"~\cite{bel2012gradual,
%meron2018patterns, 
scheffer2001catastrophic,zelnik2018regime,zelnik2015gradual}. So far, mathematical approaches to tipping have focused on examples and theory of  instability in  nonautonomous ordinary differential equation (ODE) models~\cite{alkhayuon2018rate,ashwin2017parameter, ashwin2012tipping,chen2019noise,kuehn2011mathematical,kuehn2013mathematical, o2020tipping,  ritchie2016early,ritchie2017probability,thompson2011climate,wieczorek2023rate}. These studies have identified different critical factors for tipping as well as  various tipping mechanisms. On the other hand, tipping  in spatially-extended systems modelled by nonautonomous partial differential equations (PDEs) has been much less explored~\cite{berestycki2009can,chen2015patterned,liu2014pattern}. While PDE models will likely exhibit new critical factors %(e.g. {\em critical spatial extent}) 
and interesting tipping mechanisms, their analysis is more challenging and requires new methods.

In this work, we analyze  tipping in spatially-extended systems modelled by \emph{reaction-diffusion equations} (RDEs) with reaction terms that are 
space dependent (\emph{heterogeneous}),
decaying at the boundaries ({\em bi-asymptotically  homogeneous}), and possibly time dependent (\emph{nonautonomous}).
Specifically, we develop a mathematical framework that  allows us to analyze tipping in such RDEs in terms of intersecting invariant manifolds of saddle equilibria for a suitably compactified moving-frame ODE. 
Inspired by~\cite{berestycki2009can}, 
we apply our framework to an RDE model of a geographically shifting  or shrinking ecosystem and describe two different tipping mechanisms that are characteristic of spatially-extended systems.

When discussing tipping points in nonautonomous systems with time-varying external inputs, it is useful to consider the corresponding {\em frozen system} with fixed-in-time inputs. In the frozen system, we identify a desired stable state and refer to this state as the {\em base state}.
When the external input changes over time,  the shape and position of the base state may change too, and the nonautonomous system will try to adapt to the changing base state. In other words, the nonautonomous system will try to {\em track} the stable branch of base states for the frozen system.
However, in some cases tracking is not possible, and the nonautonomous system  tips from the base state to a different
state, such as an \emph{alternative stable state}.
For example, the base state may lose stability or disappear in a classical bifurcation at some {\em critical level} of the external input. If this bifurcation is {\em dangerous}~\cite{ashwin2012tipping}, meaning that it gives rise to a discontinuity in the stable branch of base states, we say the system undergoes {\em bifurcation-induced tipping} or {\em B-tipping}~\cite{ashwin2012tipping}.
What is more, if the external input changes faster than some {\em critical rate}, the nonautonomous system may deviate too far from the changing base state, cross some threshold, and tip to an alternative stable state, even though the base state never loses stability or disappears.
Such transitions are caused entirely by the rate of change of the external input, and we say the system undergoes {\em rate-induced tipping} or {\em R-tipping}~\cite{ashwin2012tipping,o2020tipping,wieczorek2023rate}.

In our framework, we consider a one-dimensional RDE with a bi-asymptotically homogeneous reaction term.
In the moving frame,
such an RDE reduces to a special {\em nonautonomous moving-frame ODE} that is often described as a {\em bi-asymptotically autonomous} ODE~\cite{markus2016ii,wieczorek2021compactification}.
We exploit the asymptotic properties of the nonautonomous  moving-frame ODE and use the compactification technique of~\cite{wieczorek2021compactification} to
reformulate it
%the nonautonomous moving-frame ODE 
into an autonomous ODE on a suitably extended and compactified phase space. We refer to the ensuing autonomous ODE  as \emph{the compactified system}.
Different compactification methods have been used before to study the phase space near infinity \cite{dumortier2006compactification,giraldo2017saddle,giraldo2020computing,krauskopf1997bifurcations,matsue2018blow}, compute linear spectra of nonlinear wave solutions \cite{jones1990topological,kapitula2004eigenvalues}, and facilitate analysis of nonautonomous ODEs \cite{wieczorek2023rate,wieczorek2021compactification}.
A particular advantage of our compactified system is that, unlike the nonautonomous moving-frame ODE, it is autonomous and contains regular equilibria from infinity. This allows us to obtain pulse and front solutions for the original RDE by computing heteroclinic  orbits that connect an equilibrium from negative infinity to an equilibrium from positive infinity in the compactified system.
%Moreover, 
We compute such heteroclinic orbits as intersections of  invariant manifolds of these equilibria,
which automatically allows us to capture multiple coexisting pulse and front solutions. In practice,  we obtain the unstable invariant manifold of the equilibrium from negative infinity and  the stable invariant manifold of the equilibrium from positive infinity by combining an adaptive collocation method~\cite{russell1978adaptive} and pseudo-arclength continuation~\cite{doedel2007lecture}. 
Intersections of these manifolds
%, which correspond to connecting heteroclinic orbits in the compactified system, 
%{\ch \sout{which in turn correspond to pulse and front solutions in the original RDE,}} 
are detected using Lin's method for connecting orbits~\cite{lin1990using}. 
For convenience, all three numerical methods are implemented in the continuation software package AUTO~\cite{Doedel2007AUTO}. This implementation enables us to perform numerical continuation of heteroclinic orbits in the compactified system, giving rise to 
%one-parameter and two-parameter 
{\em bifurcation diagrams} of pulse and front solutions for the original RDE.
Finally, we identify tipping points as dangerous bifurcations of pulse and front solutions in the moving frame of reference.

Inspired by the work done in~\cite{berestycki2009can}, we choose as an illustrative example for our framework a mathematical model of a habitat that is shrinking in size due to, for example, increased human activity, or shifting in space due to, for example, global warming.
The difference is that we consider a more general reaction term.
%
%Our modelling approach is inspired by the work by done in~\cite{berestycki2009can} but we use a different approach for the reaction term..
%
The reaction term used in~\cite{berestycki2009can}  was a discontinuous  piecewise-homogeneous function of space. In other words, the spatial domain  was  separated into a homogeneous ``good habitat," where any non-zero population tends toward the carrying capacity, and a homogeneous ``bad habitat," where population always declines to extinction. This gives rise to a piecewise-autonomous ODE in a moving frame. Then, pulse solutions for the original RDE  were  constructed in the moving frame by gluing  orbit  segments of a linear autonomous  ODE obtained inside the bad habitat and  orbit segments of a nonlinear autonomous ODE obtained inside the good habitat.
In contrast,
our reaction term is a nonhomogeneous $C^1$-smooth function of space with a continuous transition between the good and bad habitats. Such a reaction term gives rise to a nonautonomous ODE in a moving frame, meaning that the approach of gluing orbit segments of different autonomous ODEs does not apply.
To address this problem, we propose a framework that combines compactification, invariant manifold computations, Lin's method, and numerical continuation to study tipping from pulse and front solutions in RDEs with such reaction terms.

 Most studies of geographically shifting ecosystems~ \cite{berestycki2009can,bouhours2019spreading, li2014persistence,maciel2013individual,pease1989model} focused on a monostable logistic population growth model inside the good habitat.
 On the other hand, some studies \cite{harsch2017moving,roques2008population} considered a bistable growth model that takes into account the effect of \emph{undercrowding} at low population density, also known as the Allee effect (see~\cite{volterra1938population} and~\cite[Sec.3]{dennis1989allee}).
Roques et al.~\cite{roques2008population} considered three different configurations of a two-dimensional spatial domain and found that the populations subject to the Allee effect are more sensitive to the shape and position of the habitat.
Harsch et al.~\cite{harsch2017moving} used integrodifference equations to conduct case studies for the impact of moving habitats on (i) populations subject to the Allee model, (ii) interspecific competitions, and (iii) disease-infected populations.
Other work \cite{berestycki2014can,potapov2004climate} focused on interspecific competitions and investigated the effect of moving habitats in \emph{invasion problems}.
Here, we introduce a non-homogeneous $C^1$-smooth habitat function, couple it with the Allee growth model, and 
highlight the key differences from the the logistic growth model.

As the main result, we uncover B-tipping to extinction below a \emph{critical length} of the habitat and R-tipping to extinction above a \emph{critical speed} of the shifting habitat.
Each tipping point corresponds to a fold, or equivalently saddle-node, bifurcation of pulse solutions for the RDE, which are obtained  by computing a codimension-one heteroclinic orbit along a (quadratic) tangency of invariant manifolds in the compactified system. 
Furthermore, we continue these heteroclinic orbits in the system parameters to produce two-parameter {\em tipping diagrams}, revealing nonobvious dependence of the critical length and critical speed on the diffusion, or equivalently the dispersal rate, of the population.

The organization of the paper is as follows.
\Cref{sec:model} presents the general form of the RDE and the specifics of the habitat model, and demonstrates the presence of both B-tipping and R-tipping via direct numerical simulations.
In \Cref{sec:GenApproach}, we outline the details of our mathematical framework for studying pulse and front solutions in bi-asymptotically homogeneous RDEs.
The framework is presented in four steps: (i) nondimensionalization, (ii) reformulation of the problem and comparison to practiced approaches to obtaining pulse and front solutions by computing connecting orbits in the moving-frame ODE, (iii) compactification, and (iv) numerical implementation.
In \Cref{sec:mhresults} , we demonstrate the results for pulse solutions and their bifurcations in the habitat model. 
Conclusions and final remarks are discussed in \Cref{sec:discussion}.

\section{The model}
\label{sec:model}
This paper considers nonlinear dynamics of 
RDEs in one spatial dimension
\begin{equation} 
\label{eq:main}
    u_t = D \, u_{xx} + f\left(u,H(x-ct)\right),
\end{equation}
with Dirichlet boundary conditions on an unbounded domain:

\begin{equation} 
\label{eq:DirichletBCs}
    \displaystyle{\lim_{x\to \pm \infty}} u(x,t) = u^\pm \in\mathbb{R},
\end{equation}
where the independent variables $x \in \mathbb{R}$ and $t \in [0,\infty)$ represent space and time, respectively, and the subscripts represent partial derivatives $u_t=\partial u/\partial t$ and $u_{xx} = \partial^2 u/\partial x^2$.

As an illustrative example, we consider a conceptual model of a habitat that can shrink or shift~\cite{berestycki2009can},
where the state variable  $u(x,t)\ge 0$
represents the spatiotemporal density of the inhabiting population. In this model, the constant diffusion coefficient $D\ge 0$ quantifies the magnitude of population flux from higher to lower density areas, while
the space-dependent (heterogeneous) and time-dependent (nonautonomous) 
{\em reaction term}
\begin{equation}
\label{eq:rt}
f(u,H(x-ct)),
\end{equation} 
describes population growth. The spatial extent of the habitat supporting population growth, and the linear shift of the habitat in $x$ at a given constant speed $c \ge 0$, are specified by the {\em habitat function}
$H(x-ct)$, which is introduced in the next section.
The model variables and different parameter values, along with their physical units are summarized in~\cref{table:units}.

In typical reaction-diffusion problems, the reaction term is homogeneous (space independent)
and autonomous (time independent), the boundary conditions specify what types of  traveling-wave solutions (e.g., fronts, pulses or wave trains) are possible, and the primary aim is to obtain such solutions and determine their unknown speed. 
In our problem, the boundary conditions~\eqref{eq:DirichletBCs} 
together with the reaction term~\eqref{eq:rt} specify either travelling-pulse  ($u^-=u^+$) or travelling-front ($u^-\ne u^+$)  solutions and, in contrast to the typical problems,  the speed of travelling pulses is already given by $c$.
%\footnote{{\sw(SW: need a reference here.)} {\ch I thought about this footenote and I think it needs to be removed. We can talk about it in the next meeting.} The speed of travelling fronts is given by an ``invasion speed"~\cite{}  $c^*$ that needs to be determined if $c>c^*$, or by $c$ if $c < c^*$.}
Thus, it is convenient to define the {\em moving-frame coordinate}
$$
\xi=x-ct,
$$
together with an equivalent state variable
$$
U(\xi,t) = u(x,t),
$$
and reformulate the BVP~\eqref{eq:main} and~\eqref{eq:DirichletBCs} in the moving frame with the given speed $c$, in terms of $t$ and $\xi$ as the new independent variables.
This gives  the advection-reaction-diffusion equation (ARDE)
\begin{align}
\label{eq:mainmf1}
    U_t = D \, U_{\xi \xi} + c\, U_{\xi} + f(U,H(\xi)),
\end{align}
with the boundary conditions

\begin{align}
\label{eq:mainmf2}
    \displaystyle{\lim_{\xi \to \pm \infty}} U(\xi,t) =  U^\pm \in\mathbb{R}.
\end{align}
Note that the reaction term $f$  in~\eqref{eq:mainmf1} is heterogeneous in $\xi$ but no longer depends on time $t$.
Furthermore, travelling-pulse (travelling-front) solutions with speed $c$ in the original-frame BVP~\eqref{eq:main}--\eqref{eq:DirichletBCs} correspond to stationary-pulse (stationary-front) solutions
in the moving-frame BVP~\eqref{eq:mainmf1}--\eqref{eq:mainmf2}. We will
use the term {\em pulse solutions} ({\em front solutions}) to mean either, depending on the context. 
Such solutions can be obtained by setting $U_t=0$ in~\eqref{eq:mainmf1} and solving the ensuing heterogeneous BVP
\begin{align}
\label{eq:ODEmoving}
    D \, U_{\xi \xi} + c\, U_{\xi} + f(U,H(\xi))=0,\\
\label{eq:BCmoving}
    \displaystyle{\lim_{\xi \to \pm \infty}} U(\xi) =  U^\pm \in\mathbb{R},
\end{align}
where the subscripts denote new ordinary derivatives $U_\xi=dU/d\xi$ and $U_{\xi \xi} = d^2U/d\xi^2$.
The second-order ODE in~\eqref{eq:ODEmoving} is often referred to as a {\em moving-frame ODE}. 
The method of computing pulse and front  solutions (i.e., solving the BVP~\eqref{eq:ODEmoving}--\eqref{eq:BCmoving})  depends on the form of the reaction term, as will be explained in Section~\ref{sec:mfode}.

\begin{table}[t]
\centering
\begin{tabular}{|l|l|l|l|l|}
\hline
\textbf{Quantity} & \textbf{Allee unit}  & \textbf{Allee value} &  \textbf{Logistic unit} & \textbf{Logistic value} \\ \hline
$u$, $U$       & $indiv/km$            & $[0,\infty)$ & $indiv/km$ & $[0,\infty)$  \\ \hline
$x$       & $km$            & $(-\infty,\infty)$ &  $km$ &  $(-\infty,\infty)$ \\ \hline
$t$       & $yr$            & $[0,\infty)$    & $yr$ & $[0,\infty)$ \\ \hline
$\xi$       & $km$            & $(-\infty,\infty)$ &  $km$ &  $(-\infty,\infty)$ \\ \hline
$D$       & $km^2/yr$ $km^3/yr$   & varied & $km^2/yr$ & varied \\ \hline
$\beta$   & $1/yr$      & 0.45 & $1/yr$  & 5 \\ \hline
$\gamma$  &  $km^2/(indiv^2 \ yr)$     & 1 & $km/(indiv \ yr)$ & 1 \\ \hline
$\lambda$  & $km/(indiv \ yr)$        & $4\sqrt{\beta \gamma}$ & $1/yr$ & $2\beta$ \\ \hline
$L$       & $km$            & varied & $km$ & varied\\ \hline
$a$       & $km^{-1}$       & 5 & $km^{-1}$ & 5 \\ \hline
$c$       & $km/yr$    & varied & $km/yr$ & varied \\ \hline
\end{tabular}
\caption{
Physical quantities for system \eqref{eq:main} with the logistic reaction term \eqref{eq:LogisticFunction} and the Allee reaction term \eqref{eq:AlleeFunction}; $indiv$, $km$ and $yr$ denote individuals, kilometers, and years, respectively.
}
\label{table:units}
\end{table}

\subsection{ The habitat model}
%\subsection{Heterogeneous reaction terms}
In the illustrative example of a changing habitat, we consider two distinct
population growth models that are characterized by two different reaction terms.
The logistic growth model, which is the focus of~\cite{berestycki2009can},  is characterized by the {\em logistic reaction term}
\begin{equation} \label{eq:LogisticFunction}
    f_L(U,H(\xi))=  -\beta\, U + \lambda\, H(\xi)\,  U - \gamma\, U^2,
\end{equation}
which accounts for limited resources at large population density. This quadratic function has two roots. The zero root corresponds to extinction, while the positive root corresponds to the carrying capacity of the habitat.
In contrast, the Allee growth model is characterized by the {\em Allee reaction  term}
\begin{equation} 
\label{eq:AlleeFunction}
    f_A(U,H(\xi))= -\beta \, U + \lambda\, H(\xi)\, U^2 - \gamma\, U^3, 
\end{equation}
which accounts for limited resources at large population density as well as for the undercrowding Allee effect at low population density; see~\cite{volterra1938population} and~\cite[Sec.3]{dennis1989allee}. The main difference from the logistic growth model is that  this cubic function has three roots and may give rise to bistability between extinction and carrying capacity, which are separated by the unstable Allee threshold for population growth.
Our focus is on the analysis of the Allee growth model and how it contrasts with the logistic growth model. The logistic growth model is discussed in the appendix.

\begin{figure}[t]
    \centering
    \includegraphics[scale=1.0]{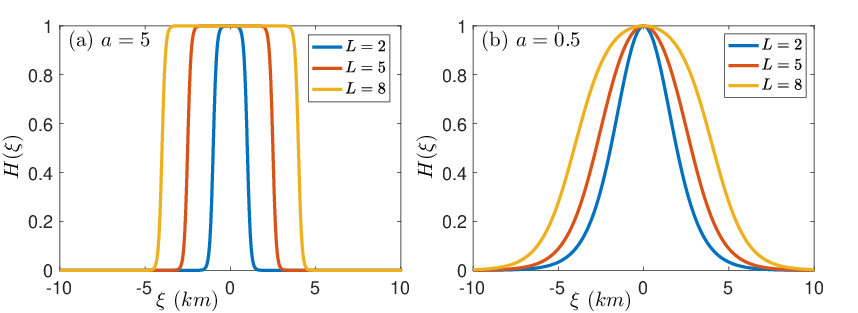}
    \caption{
    Continuous heterogeneous-in-$\xi$ habitat function $H\big( \xi)$ as defined in~\eqref{eq:habitat} with (a) a steep spatial gradient $a=5\ km^{-1}$ and (b) a gentle spatial gradient $a=0.5 \ km^{-1}$. 
    In each panel, the three different colours represent three different lengths $L=2,5$, and $8 \ km$.
    } 
    \label{fig:habitat}
\end{figure}

For both growth models, the constant parameters  $\beta  \ge 0$ and $\gamma  \ge 0$ represent the linear and nonlinear death rates, respectively.
The second term  in~\eqref{eq:LogisticFunction} and~\eqref{eq:AlleeFunction} characterises birth processes and consists of two factors. The constant parameter $\lambda \ge 0$
is the birth rate\,\footnote{Note that $\lambda$ is the linear birth rate in the logistic model~\eqref{eq:LogisticFunction}, and nonlinear birth rate in the Allee model~\eqref{eq:AlleeFunction}.}, and the dimensionless and  heterogeneous-in-$\xi$ habitat function
\begin{equation}
\label{eq:habitat}
    H(\xi) = \dfrac{\tanh\left(a\left(\xi+ L/2\right) \right)-\tanh\left(a\left(\xi-L/2\right)\right)}{2\tanh\left(aL/2\right)},
\end{equation}
specifies the position of the {\em good habitat} patch in the moving frame; see~\Cref{fig:habitat}.
By good habitat, we mean the $\xi$-interval for which $H \approx 1$. We also use the terms {\em bad habitat} to refer to the two $\xi$-intervals for which $H \approx 0$ and {\em transitional habitat} to refer to the two $\xi$-intervals where $0\lesssim H \lesssim 1$.
Here, $a>0$ quantifies the spatial slope of the transitional habitat, and $L>0$ approximates the length of the good habitat when  $a$ is large enough.

%{\sw (Suggestion: Make this ``2.2 Rate-induced tipping of a shifting habitat", start with a one-parameter bifurcation diagram from Ruaidhri's thesis, then give the motivating discussion below.}
To motivate our work, we consider two cases for a given $L$.
In the case of sufficiently large $a$,
there is an abrupt transition between the good and bad habitats, with 
the length of the good habitat  being approximately $L$, and with a relatively short length of the transitional habitat; see \Cref{fig:habitat}(a). 
In this case, the heterogeneous habitat function in~\eqref{eq:habitat}
can be approximated by a piecewise-homogeneous function which greatly simplifies analysis of pulse solutions -- we explain this in more detail in Section~\ref{sec:mfode}.
However, in the case of sufficiently small $a$, the transitional habitat extends over relatively wide $\xi$-intervals, and the length of the good habitat is noticeably shorter than $L$; see \Cref{fig:habitat}(b). This means that the piecewise-homogeneous approximation is no longer valid, and there is a need for an alternative approach to analyze pulse solutions.

\subsection{B-tipping and R-tipping in the habitat model}
\label{sec:tippingInto}
\begin{figure}[t]
   \centering
    \includegraphics[scale=1.0]{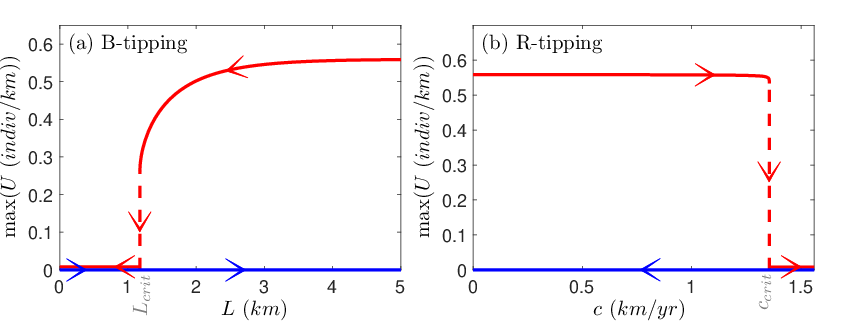}
    \caption{
    One-parameter attractor diagrams of pulses for the RDE~\eqref{eq:main}--\eqref{eq:DirichletBCs} with the Allee reaction term~\eqref{eq:AlleeFunction}, obtained using the ARDE~\eqref{eq:mainmf1}--\eqref{eq:mainmf2} with~\eqref{eq:AlleeFunction}.
    (a) Decreasing (red) and then increasing (blue) the length $L$ of a static habitat ($c=0 \ km/yr$) shows B-tipping  from a standing pulse to extinction at a critical length $L\approx 1.1748727 \ km$.
    (b) Increasing (red) and then decreasing (blue) the speed $c$ of a moving habitat with a fixed $L_{crit}=5 \ km$ shows R-tipping from a travelling pulse to extinction at a critical speed $c_{crit}\approx 1.35536596 \ km/yr$. Other parameter values are given in \Cref{table:units}.
    }
    \label{fig:BtippingAndRtipping}
\end{figure}

To give a taste of different tipping mechanisms that are present in the habitat model, we perform direct numerical simulations\footnote{Using the method of lines~\cite{hamdi2007method,zafarullah1970application}.} to detect stable pulse solutions in the ARDE~\eqref{eq:mainmf1}--\eqref{eq:mainmf2} with the Allee reaction term~\eqref{eq:AlleeFunction} and $U^{\pm}=0$; see~\Cref{fig:BtippingAndRtipping}.

In \Cref{fig:BtippingAndRtipping}(a), we set $c=0$ and simulate a slowly shrinking habitat. We start with $L=5$, detect a stable standing pulse solution, then decrement $L$, use the previously obtained solution as an initial condition, converge to the new stable solution, and repeat this procedure until we reach $L=0$. The result is a stable branch of standing pulses that terminates at
the {\em critical length} $L=L_{crit} \approx 1.17537$. For $0< L < L_{crit}$, the system always converges to the extinction solution $U(x)=0$. Subsequently, we increment $L$ using the same procedure, which shows that the branch of extinction solutions is stable, at least up to $L=5$. Thus, the shrinking habitat undergoes B-tipping to extinction upon decreasing $L$, and this transition  cannot be reversed by increasing $L$ back to its initial value.

In \Cref{fig:BtippingAndRtipping}(b), we fix $L=5$ and simulate a geographically shifting habitat at different speeds $c>0$. We start with $c=0$, detect a stable standing pulse solution, then increment $c$, and proceed in the same manner as in \Cref{fig:BtippingAndRtipping}(a). The result is a stable branch of travelling pulses that terminates at
the {\em critical speed} $c=c_{crit} \approx 1.35527$. For $c > c_{crit} $, the system always converges to the extinction solution $U(x)=0$. Thus, the moving habitat undergoes R-tipping to extinction above the critical rate
$c_{crit}$. Note that, in contrast to the branch of standing pulses near $L_{crit}$ in (a), the branch of travelling pulses in (b) remains nearly constant, showing no indication of the imminent critical speed.
%\todo{Early warning is used in the sense of ``time before tipping" and does not apply to this figure as $c$ is not increased over time}
%
In the remainder of the paper, we develop a framework to study pulse solutions in bi-asymptotically homogeneous RDEs and use this framework to uncover and discuss the dynamical mechanisms responsible for both tipping examples in \Cref{fig:BtippingAndRtipping}.

\section{The framework} \label{sec:GenApproach}
We here propose a framework that facilitates analysis of travelling pulses and fronts in the heterogeneous and nonautonomous
RDE~\eqref{eq:main} or stationary fronts and pulses in the heterogeneous ARDE~\eqref{eq:mainmf1}.
This framework
\begin{itemize}
    \item 
    is applicable to $C^1$-smooth
    reaction terms that are bi-asymptotically homogeneous in $\xi$, meaning that\,\footnote{Note that this applies to  the given habitat function in~\eqref{eq:habitat} with any combination of $a$ and $L$.}
    \begin{equation*}
     H(\xi) \to h^\pm\in\mathbb{R}\;\;\;\mbox{and}\;\;\;
    f(U,H(\xi))\to f(U,h^\pm)\;\;\;\mbox{as}\;\;\;\xi\to\pm\infty,
    \end{equation*}
    and Dirichlet boundary conditions 
    $$
    \lim_{\xi\to \pm \infty}U(\xi,t) = U^\pm \in\mathbb{R};
    $$
     \item 
     uses the compactification technique of~\cite{wieczorek2021compactification} to 
     transform a pulse or front solution in system~\eqref{eq:main} or~\eqref{eq:mainmf1} into a heteroclinic orbit in a suitably compactified system;
     \item
     uses Lin's method~\cite{lin1990using} implemented in the continuation software package AUTO~\cite{Doedel2007AUTO,krauskopf2008lin} to compute heteroclinic orbits in the compactified system that correspond
to pulse and front solutions in system~\eqref{eq:main} or~\eqref{eq:mainmf1}.
\end{itemize}
We introduce this framework in four steps.

\subsection{Nondimensionalization}

\begin{table}[t]
\centering
\begin{tabular}{|l|l|l|l|l|}
\hline
\textbf{Quantity} & \textbf{Allee rescaling}  & \textbf{Allee value} &  \textbf{Logistic rescaling} & \textbf{Logistic value} \\ \hline
 $\mathcal{U}$        & $\sqrt{\frac{\gamma}{a^2 D}} \, U$           & $[0,\infty)$ & $\frac{\gamma}{a^2 D} \, U$  & $[0,\infty)$  \\ \hline
$z$       & $a \, \xi$            & $(-\infty,\infty)$ &  $a \, \xi$ &  $(-\infty,\infty)$ \\ \hline
$\tau$       & $a^2 D t$            & $[0,\infty)$    & $a^2 D t$ & $[0,\infty)$ \\ \hline 
$\tilde{\beta}$   & $\sqrt{\frac{\beta}{a^2 D}}$       & varied  & $\sqrt{\frac{\beta}{a^2 D}}$ & varied \\ \hline
$\tilde{L}$       & $a \, L$            & varied  & $a \, L$ & varied \\ \hline
$\tilde{c}$   & $\frac{c}{a D}$     & varied  & $\frac{c}{a D}$ & varied \\ \hline
\end{tabular}
\caption{\label{table:dimensionless}
Dimensionless variables and parameters for system \eqref{eq:rescaledODEmoving}--\eqref{eq:rescaledODEmoving2} with $D>0$.
}
\end{table}

In the first step, summarized in \Cref{table:dimensionless}, we rewrite the BVP~\eqref{eq:ODEmoving}--\eqref{eq:BCmoving} in terms of a dimensionless moving-frame coordinate $z$ and a dimensionless state variable $\mathcal{U}(z)$ as\,\footnote{Note the slight abuse of notation where we use the same symbol $\mathcal{U}(z)$ for differently rescaled $U(\xi)$ in the logistic and Allee models. 
}
\begin{align}
     \mathcal{U}_{z z} + \tilde{c} \, \mathcal{U}_{z} + \tilde{f}(\mathcal{U},\tilde{H} (z) )=0, \label{eq:rescaledODEmoving}\\
    \displaystyle{\lim_{z\to \pm \infty}} \mathcal{U}(z) = {\mathcal{U}}^\pm \in\mathbb{R}.
    \label{eq:rescaledODEmoving2}
\end{align}
A particular advantage of this nondimensionalization is that 
the number of parameters in the system reduces from seven to just three, namely, $\tilde{c}$, $\tilde{\beta}$, and $\tilde{L}$.
%\sout{Two particular advantages of this nondimensionalization are: (i) the number of parameters in the system reduces from seven to just three, namely, $\tilde{c}$, $\tilde{\beta}$ and $\tilde{L}$, and (ii) stability analysis of equilibria in the compactified system in \cref{sec:stabilityOfEq} ahead is simplified.}
This advantage becomes clear from the rescaled Allee reaction term~\eqref{eq:AlleeFunction},
\begin{equation}
\label{eq:fAr}
\tilde{f}_A(\mathcal{U},\tilde{H}(z)) = 
\mathcal{U}
\left(
-\tilde{\beta}^2
+ 4 \tilde{\beta} \, \tilde{H}(z)\, \mathcal{U}
- \mathcal{U}^2
\right),
\end{equation}
and the rescaled habitat function~\eqref{eq:habitat},
 \begin{equation}
 \label{eq:habitatresc}
    \tilde{H}(z) = \dfrac{\tanh(z+\tilde{L}/2 )-\tanh(z-\tilde{L}/2)}{2\tanh(\tilde{L}/2)}.
\end{equation}

\subsection{Pulses and fronts in the moving-frame ODE}
\label{sec:mfode}

\begin{figure}[t]
\label{fig:schematic}
    \centering
    \includegraphics[scale=0.8]{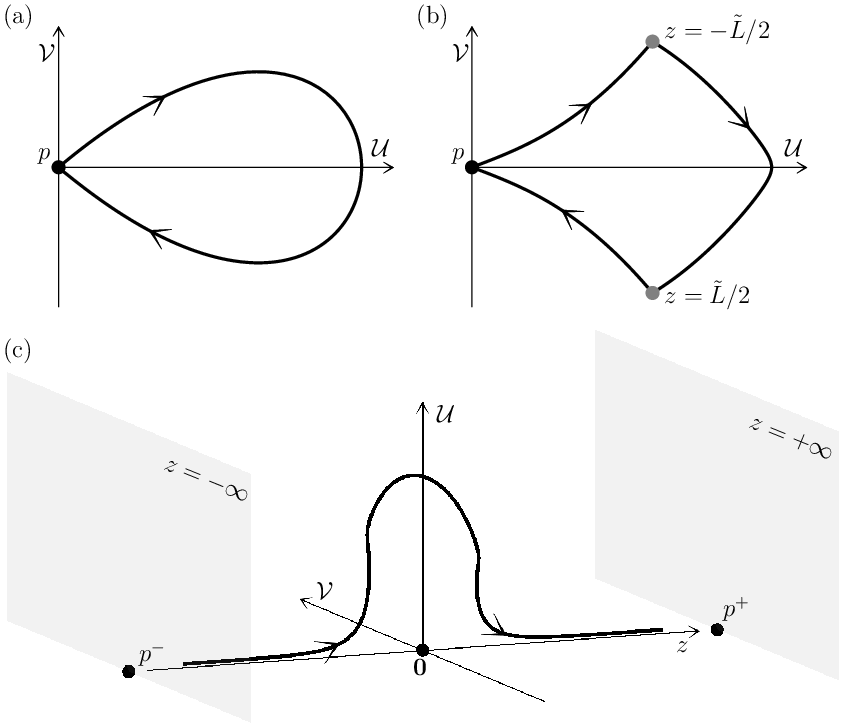}
    \caption{
    A pulse solution for an RDE  can be obtained in the moving-frame ODE by computing (a)~a smooth homoclinic orbit to a saddle $p$ when the reaction term is homogeneous ($z$-independent), (b) a piecewise-smooth homoclinic orbit to a saddle $p$ when the reaction term is piecewise-homogeneous (piecewise-constant in $z$), and (c) a trajectory that limits to an equilibrium $p^\pm$ as $z$ tends to $\pm\infty$ when the reaction term is bi-asymptotically homogeneous.
    }
\end{figure}

In the second step, we rewrite the second-order BVP~\eqref{eq:rescaledODEmoving}--\eqref{eq:rescaledODEmoving2} as a first-order BVP at the expense of introducing an additional dependent variable $\mathcal{V}(z)={\mathcal{U}}_z(z)$:
\begin{align}
           \mathcal{U}_z & =  \mathcal{V}, 
           \label{eq:1stOrderU} \\
           \mathcal{V}_z & =  - \tilde{c} \, \mathcal{V} - \tilde{f} (\mathcal{U},\tilde{H}(z) ), 
           \label{eq:1stOrderV}\\
           \displaystyle{\lim_{z\to \pm \infty}} \mathcal{U}(z) & = {\mathcal{U}}^\pm \in\mathbb{R}.
           \label{eq:1stOrderBC}
\end{align}
We then view the moving-frame ODEs~\eqref{eq:1stOrderU} and~\eqref{eq:1stOrderV} as a dynamical system on $\mathbb{R}^2$, where $z$ plays the role of time. 
The following three paragraphs overview how one can obtain pulse and front
solutions to the BVP~\eqref{eq:1stOrderU}--\eqref{eq:1stOrderBC} depending on the nature of the reaction term $\tilde{f}$.

In typical problems with a $z$-independent reaction term  $\tilde{f}(\mathcal{U}, h)$  and a constant $h \in\mathbb{R}$,
the ensuing dynamical system is \emph{autonomous}.
Pulse solutions are possible if 
${\mathcal{U}}^\pm = {\mathcal{U}}^*$ and there is a saddle equilibrium point $p=({\mathcal{U}}^*,0)$ in the $(\mathcal{U},\mathcal{V})$ phase plane,
in which case such solutions can be computed as homoclinic orbits to $p$; see~\Cref{fig:schematic}(a).
Similarly, front solutions are possible if ${\mathcal{U}}^-\ne {\mathcal{U}}^+$ and there are two different  equilibrium points, 
 $p^- = ({\mathcal{U}}^-,0)$ and $p^+ = ({\mathcal{U}}^+,0)$,
in which case such solutions can be computed as heteroclinic orbits
from ${p}^-$ to ${p}^+$. This typical approach has been widely implemented, for example, in~\cite[Ch.7]{kuehn2019pde}, \cite[Ex.6.3]{kuznetsov2013elements} and~\cite{champneys1998homoclinic,jones1984stability,sandstede2002stability,sandstede1999essential}.

In problems where  the reaction term  $\tilde{f}(U,\tilde{H}(z))$ can be approximated by a piecewise-homogeneous function,
the ensuing  dynamical system is \emph{piecewise-autonomous}. For example, if
%(see fig.~\ref{fig:habitat}(a)),
\begin{equation}
\label{eq:habitatPWC}
   \tilde{f}(U,\tilde{H}(z)) =
  \begin{cases}
    \tilde{f}(U, h^-), \quad h^- \in \mathbb{R}, & \text{if}\;\;\; z < -\tilde{L}/2 , \\
    \tilde{f}(U, h^0), \quad \ h^0 \in \mathbb{R},  & \text{if } -\tilde{L}/2\le z \le \tilde{L}/2, \\
    \tilde{f}(U, h^+), \quad h^+ \in \mathbb{R},  & \text{if}\;\;\; z > \tilde{L}/2,
  \end{cases}
\end{equation}
the system is given by a set of three  (constituent) autonomous dynamical systems that are defined on three adjacent $z$-intervals with shared boundaries at $z=\pm L/2$; see~\Cref{fig:habitat}(a).
 Pulse solutions are possible if  ${\mathcal{U}}^\pm=\mathcal{U}^*$ and there is a saddle equilibrium point $p=(\mathcal{U}^*,0)$ for $z<-\tilde{L}/2$ and $z>\tilde{L}/2$. Then, such solutions can be computed as {\em piecewise-smooth homoclinic orbits} to  $p$, via a concatenation of three orbit segments of the three constituent autonomous systems that match at `times' $z=\pm \tilde{L}/2$; see~\Cref{fig:schematic}(b).
Similarly, front solutions are possible if ${\mathcal{U}}^-\ne {\mathcal{U}}^+$ and there are two different equilibrium points,
$p^- = ({\mathcal{U}}^-,0)$
for $z<-\tilde{L}/2$ and $p^+ = ({\mathcal{U}}^+,0)$  for $z>\tilde{L}/2$, in which case such solutions can be computed as {\em piecewise-smooth heteroclinic orbits} from $p^-$ to $p^+$.
Such a concatenation technique  was proposed  in~\cite{jones1986instability,jones1988instability,langbein1985generalized} and has been implemented in various contexts and applications~\cite{berestycki2009can,kozak2019pattern,nishi2013dynamics,van2010pinned}.

By contrast, a general $z$-dependent reaction term $\tilde{f}(\mathcal{U},\tilde{H}(z))$ poses an obstacle to computing pulse and front solutions: it gives rise to a {\em nonautonomous} dynamical system~\eqref{eq:1stOrderU}--\eqref{eq:1stOrderV} that has no equilibrium points in the extended $(\mathcal{U},\mathcal{V},z)$ phase space;
%\,\footnote{Note that the stationary solution $p(t)=(\mathcal{U}(t),\mathcal{V}(t)) =(0,0)$ becomes a line in the extended $(\mathcal{U},\mathcal{V},s)$ phase space.}; 
see \Cref{fig:schematic}(c).  This obstacle becomes particularly apparent when the BVP~\eqref{eq:1stOrderU}--\eqref{eq:1stOrderBC} has multiple `nearby' pulse or front solutions
that may be difficult to capture by a shooting method or a collocation method.
%Thus requires an alternative approach.
To overcome this obstacle, we exploit the fact that
the reaction term is bi-asymptotically homogeneous in the sense that
$$
\tilde{f} (\mathcal{U},\tilde{H}(z) )\to \tilde{f} (\mathcal{U},h^\pm )\quad\mbox{as}\quad z\to\pm\infty.
$$
Specifically, we use an equilibrium point  
\begin{equation}
\label{eq:p-}
 p^- =  (\mathcal{U},\mathcal{V}) = ({\mathcal{U}}^-,0),\nonumber
\end{equation}
for the autonomous past limit system
\begin{align}
\begin{split}
           \mathcal{U}_z  &=  \mathcal{V},\\
           \mathcal{V}_z & =  - \tilde{c} \, \mathcal{V} - \tilde{f} (\mathcal{U},h^-),
\end{split}
\label{eq:pastSystem}
\end{align}
and an equilibrium point  
\begin{equation}
\label{eq:p+}
 p^+ =  (\mathcal{U},\mathcal{V}) = ({\mathcal{U}}^+,0),\nonumber
\end{equation}
for the autonomous future limit system
\begin{align}
\begin{split}
           \mathcal{U}_z  &=  \mathcal{V},\\
           \mathcal{V}_z & =  - \tilde{c} \, \mathcal{V} - \tilde{f} (\mathcal{U},h^+ ),
\end{split}
\label{eq:futureSystem}
\end{align}
to construct pulse or front solutions as  heteroclinic orbits from
 $p^-$ to $p^+$ in a suitably compactified system.

\subsection{Compactification} 
\label{sec:compact}

In the third step, we bring in equilibria of the limit systems from infinity by reformulating the nonautonomous system \eqref{eq:1stOrderU}--\eqref{eq:1stOrderV} on $\mathbb{R}^2$ 
into an autonomous compactified system on $\mathbb{R}^2\times [-1,1]$.
This requires a suitable coordinate transformation that makes the additional dependent variable bounded and ensures that the compactified system is at least $C^1$-smooth on the extended phase space.

Reference~\cite[Sec. 4]{wieczorek2021compactification} constructs examples of coordinate transformations for different asymptotic decays of the  nonautonomous reaction term, ranging from sub-logarithmic to super-exponential decays. 
Here, we  focus on the case where the nonautonomous reaction term decays exponentially with a {\em decay coefficient} $\rho>0$ as $z$ tends to $\pm\infty$, in the sense that~\cite{wieczorek2023rate}
\begin{equation*}
  \lim_{z\to \pm\infty}\frac{\tilde{H}_z(z)}{e^{\mp\rho z}}\;\;\mbox{exists for some $\rho>0$}.
  % for some}\;\;\rho>0, 
  %\label{eq:smoothnessCon}
\end{equation*}
Thus, we use the parameterized coordinate transformation~\cite[Eq.(48)]{wieczorek2021compactification}, designed for exponentially or faster decaying reaction terms, to augment the nonautonomous system~\eqref{eq:1stOrderU}--\eqref{eq:1stOrderV} with
\begin{equation}
    s = g_{\alpha}(z) = \tanh\left(\frac{\alpha}{2} z\right),
    \label{eq:compactificationFunction}
\end{equation}
as a third dependent variable.\footnote{ Note that nonautonomous reaction terms with algebraic or logarithmic decay will require different transformations. Also note that we use the subscript $\alpha$ to denote dependence on $\alpha$, not a partial derivative.}  
The {\em compactification parameter} $\alpha>0$ quantifies the rate of the exponential decay of both $s(z)$ and $s_z(z)$. 
We note that $s \in (-1,1)$ for $z \in \mathbb{R}$, use the definition of $\tanh^{-1}$ to obtain the inverse transformation
\begin{equation}
\label{eq:inverseCF}
 z= 
g_{\alpha}^{-1}(s) = \frac{1}{\alpha} \ln{\frac{1+s}{1-s}},
\end{equation}
and  differentiate $g_{\alpha}(z)$ in~\eqref{eq:compactificationFunction} with respect to $z$ to derive the ODE for $s$:
$$
s_z = \frac{\alpha}{2}\,(1-s^2).
$$
Next, we continuously extend the new dependent variable $s$ to include the limits from $z=\pm\infty$, which correspond to $s=\pm 1$. 
This gives the {\em autonomous compactified system}
\begin{equation}
  \left\{
    \begin{array}{rrl}
           \mathcal{U}_z & = & \mathcal{V}, \\
           \mathcal{V}_z & = & - \tilde{c} \, \mathcal{V} - \tilde{f} \Big(\mathcal{U},\tilde{H}_\alpha(s)\Big), \\
           s_z & = & \dfrac{\alpha}{2} (1-s^2),
    \end{array} \right.
    \label{eq:compactifiedODE}
\end{equation}
defined on $\mathbb{R}^2\times[-1,1]$, with the continuously-extended habitat function
\begin{equation}
\label{eq:Hcompact}
    \tilde{H}_\alpha(s) =
  \begin{cases}
    \tilde{H}(g^{-1}_{\alpha}(s)), & \text{for}\;\;\; s \in (-1,1),  \\
    h^-, & \text{for } s=-1, \\
    h^+, & \text{for}\;\;\; s=1.
  \end{cases}
\end{equation}
%and the inverse coordinate transformation given by
%\begin{equation}
%\label{eq:inverseCF}
% z= 
%g_{\alpha}^{-1}(s) = \frac{1}{\alpha} \ln{\frac{1+s}{1-s}}.
%\end{equation}
 %Furthermore, if $\tilde{H}(z)$ decays exponentially at a rate $\rho$ as $z\to\pm\infty$, in the sense that
% \begin{remark}
%\label{Rmk:smoothnessCon}
It follows from~\cite[Cor.4.1]{wieczorek2021compactification} that if $\tilde{H}(z)$ decays exponentially with a decay coefficient $\rho > 0$,
%in the sense that
%\begin{equation*}
%  \lim_{z\to \pm\infty}\frac{\tilde{H}_z(z)}{e^{\mp\rho z}}\;\;\mbox{exists for some $\rho>0$},% for some}\;\;\rho>0, 
  %\label{eq:smoothnessCon}
%\end{equation*}
then the compactified system~\eqref{eq:compactifiedODE} is continuously differentiable on the extended phase space $\mathbb{R}^2\times[-1,1]$ for any
$$
\alpha\in(0,\rho].
$$
In other words, one needs to ensure that the compactification parameter $\alpha$ does not exceed the exponential decay coefficient $\rho$.
%\end{remark}

A particular advantage of compactification is that the flow-invariant planes of the compactified system~\eqref{eq:compactifiedODE}, namely,
%$\{s=-1\}$ and $\{s=1\}$ 
$$
S^-=\mathbb{R}^2\times\{-1\}\quad\mbox{and}\quad S^+=\mathbb{R}^2\times \{1\},
$$
contain equilibria $p^-$ and $p^+$
of the autonomous past~\eqref{eq:pastSystem} and future~\eqref{eq:futureSystem} limit systems, respectively. 
When embedded in the extended phase space of the compactified system~\eqref{eq:compactifiedODE},
%~\eqref{eq:compactifiedODE},
$p^-$ becomes 
$$
\tilde{p}^-=(\mathcal{U},\mathcal{V},s)=(\mathcal{U}^-,0,-1)\in S^-, 
$$
and gains one additional eigendirection $v^-$ with positive eigenvalue $\alpha>0$, and  $p^+$ becomes
$$
\tilde{p}^+=(\mathcal{U},\mathcal{V},s)=(\mathcal{U}^+,0,1)\in S^+,
$$
and gains one additional eigendirection $v^+$ with negative eigenvalue $-\alpha<0$; see~\cite[Rem. 3.1 and Cor. 4.1]{wieczorek2021compactification}.
Thus, a pulse or front solution to the BVP~\eqref{eq:ODEmoving}--\eqref{eq:BCmoving} can be computed as a heteroclinic connecting orbit from $\tilde{p}^-$ to $\tilde{p}^+$ in the compactified system~\eqref{eq:compactifiedODE}.
The computation of such connecting orbits becomes more convenient if
\begin{itemize}
\item [(i)]
     $v^-$  is normal to  $S^-$ and typical trajectories leave  $\tilde{p}^-$  along $v^-$,
    \item [(ii)]
     $v^+$ is normal to  $S^+$ and typical trajectories approach $\tilde{p}^+$ along $v^+$.
%    \item [(i)]
%     $v^-$  is normal to  $S^-$, $v^+$ is normal to  $S^+$, and
%    \item [(ii)]
%    typical trajectories leave  $\tilde{p}^-$  along $v^-$  and approach $\tilde{p}^+$ along $v^+$.
\end{itemize}
Suppose that ${p}^-$ and ${p}^+$ are hyperbolic and note from the discussion of additional eigenvalues due to compactification that $\tilde{p}^-$ and $\tilde{p}^+$ must be hyperbolic too.
If the unstable (stable) invariant manifold of $\tilde{p}^-$ ($\tilde{p}^+$) is of dimension one, conditions (i) and (ii) are satisfied  for any $\alpha\in(0,\rho)$; see~\cite[Rem. 3.1 and Cor. 4.1]{wieczorek2021compactification} and~\cite[Prop.6.3]{wieczorek2023rate}. 
In contrast, if the unstable invariant manifold of $\tilde{p}^-$ is of dimension greater than one, condition (i) is satisfied  for any $\alpha\in\left(0,\min\{\rho,l^-\}\right)$, where $l^-$  is the smallest-magnitude eigenvalue within the unstable eigenspace of $p^-$  in the autonomous past limit system~\eqref{eq:pastSystem}; see~\cite[Prop.6.3]{wieczorek2023rate}.
Similarly, if the stable invariant manifold of $\tilde{p}^+$ is of dimension greater than one, condition (ii) is satisfied  for any  $\alpha\in\left(0,\min\{\rho,|l^+| \}\right)$, where $l^+$ is the smallest-magnitude eigenvalue within the stable eigenspace of $p^+$ in the autonomous future limit system~\eqref{eq:futureSystem}.
%In contrast, if the unstable (stable) invariant manifold of $\tilde{p}^-$ ($\tilde{p}^+$) is of dimension greater than one, conditions (i) and (ii) are satisfied  for any $\alpha\in\left(0,\min\{\rho,l^-\}\right)$ $(\alpha\in\left(0,\min\{\rho,|l^+| \}\right))$, where $l^-$ ($l^+$) is the smallest-magnitude eigenvalue within the unstable (stable) eigenspace of $p^-$ ($p^+$) in the autonomous past~\eqref{eq:pastSystem} (future~\eqref{eq:futureSystem}) limit system; see~\cite[Prop.6.3]{wieczorek2023rate}.

\subsection{Numerical implementation} 
\label{sec:numerical}

In the fourth step,
we outline a numerical setup for computing
pulse and front solutions in a reaction-diffusion system~\eqref{eq:main} or~\eqref{eq:mainmf1} as heteroclinic orbits from $\tilde{p}^-\in S^-$ to $\tilde{p}^+\in S^+$
in the autonomous compactified system~\eqref{eq:compactifiedODE}. For convenience, we use the  continuation software package AUTO~\cite{Doedel2007AUTO}, which allows numerical continuation of solutions to autonomous ODEs subject to boundary conditions.

We start with the notation and write
$$
\Gamma(z) = \left(
\mathcal{U}(z),\mathcal{V}(z), s(z)
\right)\in\mathbb{R}^2\times[-1,1],
$$
to denote a solution to the compactified system~\eqref{eq:compactifiedODE} at `time' $z$. 
%We also re-write the compactified system~\eqref{eq:compactifiedODE} in the following vector form 
%\begin{equation}
%    \dfrac{d\Gamma}{dz} =   \mathcal{F}\left( \Gamma, \mu \right),
%    \label{eq:VectorForm}
%\end{equation}
%where $\mathcal{F}$ is the right-hand side of ~\eqref{eq:compactifiedODE} and $\mu$ is the set of system parameters.
Furthermore, we write $E^u(\tilde{p}^-)$ to denote the unstable eigenspace of  $\tilde{p}^-$, $W^u(\tilde{p}^-)$ to denote a (numerical approximation of a local) unstable manifold of $\tilde{p}^-$,  $E^s(\tilde{p}^+)$
to denote the stable eigenspace of $\tilde{p}^+$, and $W^s(\tilde{p}^+)$ to denote a (numerical approximation of a local) stable manifold of  $\tilde{p}^+$ in~\eqref{eq:compactifiedODE}.

A heteroclinic orbit from $\tilde{p}^-$ to $\tilde{p}^+$ in~\eqref{eq:compactifiedODE} is a special orbit that connects $\tilde{p}^-$ to $\tilde{p}^+$ along an intersection of $W^u(\tilde{p}^-)$ and $W^s(\tilde{p}^+)$. 
Such a connection can be
%A heteroclinic orbit from $\tilde{p}^-$ to $\tilde{p}^+$ is
 approximated by a finite-time orbit segment that starts from  $E^u(\tilde{p}^-)$ sufficiently close to $\tilde{p}^-$ at time $z=0$, and crosses $E^s(\tilde{p}^+)$ sufficiently close to $\tilde{p}^+$ at some later time $z=Z>0$.
 Specifically,  for fixed $\varepsilon^-, \varepsilon^+>0$, we consider  a finite-time orbit segment
\begin{equation}
\Gamma  := \left\{
\Gamma(z)~:~
z\in[0,Z]
\right\}\subset \mathbb{R}^2\times[-1,1], 
\label{eq:orbitSegment}
\end{equation}
where
\begin{equation}
  \left\{
    \begin{array}{ll}
           \Gamma(0) & \in E^u(\tilde{p}^-) \quad\mbox{and}\quad  
\left\|\Gamma(0) -\tilde{p}^-\right\| = \varepsilon^-, \\
           \Gamma(Z) & \in E^s(\tilde{p}^+) \quad\mbox{and}\quad 
\left\|\Gamma(Z) -\tilde{p}^+\right\| = \varepsilon^+.
    \end{array} \right.
    \label{eq:heteroBC}
\end{equation}
In other words, $\Gamma(0)$ lies on a half $(n-1)$-sphere of radius $\varepsilon^-$ about $\tilde{p}^-$  within an $n$-dimensional $E^u(\tilde{p}^-)$, and  $\Gamma(Z)$ lies on a half $(m-1)$-sphere of radius $\varepsilon^+$ about $\tilde{p}^+$ within an $m$-dimensional $E^u(\tilde{p}^-)$.\footnote{The other half of the $(n-1)$- and $(m-1)$-spheres lies outside the compactified phase space $\mathbb{R}^2\times[-1,1]$.}
Moreover, $\varepsilon^-$ and $\varepsilon^+$ are chosen small enough so that $E^u(\tilde{p}^-)$ is a good approximation of $W^u(\tilde{p}^-)$ on the half $(n-1)$-sphere  about $\tilde{p}^-$, and 
$E^s(\tilde{p}^+)$ is a good approximation of $W^s(\tilde{p}^+)$ on the half $(m-1)$-sphere of about $\tilde{p}^+$.

To obtain such finite-time orbit segments, we implement Lin's method~\cite{lin1990using} in AUTO~\cite{Doedel2007AUTO}; see also \cite{alkhayuon2018rate,giraldo2018cascades,krauskopf2008lin,mujica2017lin, musoke2020surface, oldeman2003homoclinic, zhang2012find}.
First, we need to identify a two-dimensional cross section $\Sigma$ that is transversal to the flow along the heteroclinic orbit and, for practical purposes, 
%we choose $\Sigma$ to be 
is sufficiently far from both $\tilde{p}^-$ and $\tilde{p}^+$. We note that $s_z =\alpha/2 >0$ at $s=0$ and choose
\begin{equation}
\label{eq:Sigma}
\Sigma:=\{(\mathcal{U},\mathcal{V},s)~:~ s=0 \} = \mathbb{R}^2\times \{0\},
\end{equation}
which satisfies the transversality requirement.
%Also note that our choice of $\Sigma$ is half way between $\tilde{p}^-$ and $\tilde{p}^+$,

Second, we compute two orbit segments, denoted $\Gamma^-$ and $\Gamma^+$.  We define $\Gamma^-$ as an orbit segment that starts at time $z=0$ from an $(n-1)$-sphere about $\tilde{p}^-$ within $E^u(\tilde{p}^-)$ and meets $\Sigma$
at time $z=Z^- >0$:
$$
\Gamma^- := \left\{
\Gamma^-(z)~:~
z\in[0,Z^-]
\right\}\subset \mathbb{R}^2\times[-1,0],
$$
where
\begin{align}
 \left\{
    \begin{array}{lll}
           \Gamma^-(0) & \in & E^u(\tilde{p}^-) \quad\mbox{and}\quad  
\left\|\Gamma^-(0) -\tilde{p}^-\right\| = \varepsilon^-, \\
           \Gamma^-(Z^-) & \in & \Sigma.
    \end{array} \right.
    \label{eq:leftBC}
\end{align}
Going backward in $z$, we define $\Gamma^+$ as an orbit segment that starts at `time' $z=0$ from an $(m-1)$-sphere about $\tilde{p}^+$ within $E^s(\tilde{p}^+)$ and meets $\Sigma$
at `time' $z=Z^+<0$:
$$
\Gamma^+  := \left\{
\Gamma^+(z)~:~
z\in[Z^+,0]
\right\}\subset \mathbb{R}^2\times[0,1],
$$
where
\begin{equation}
  \left\{
    \begin{array}{lll}
           \Gamma^+(0) & \in & E^s(\tilde{p}^+) \quad\mbox{and}\quad 
\left\|\Gamma^+(0) -\tilde{p}^+\right\| = \varepsilon^+, \\
           \Gamma^+(Z^+) & \in & \Sigma.
    \end{array} \right.
    \label{eq:rightBC}
\end{equation}
The computation of these two orbit segments can be performed by the shooting method, that is, by reducing the BVP to an IVP and solving the IVP using direct $z$-integration. Here, we use AUTO instead, which solves the BVP directly by combining an adaptive collocation method~\cite{russell1978adaptive} and pseudoarclength continuation~\cite{doedel2007lecture}. One advantage of our approach is to equidistribute the local discretization error along the computed orbit segment. Another advantage is that once an orbit segment satisfying~\eqref{eq:leftBC} or~\eqref{eq:rightBC} is computed in AUTO, it can be readily numerically continued in AUTO by varying a system parameter or a  {\em parameterized boundary condition.}
To be more precise, in the case where the eigenspace $E^{u}(\tilde{p}^-)$  is one-dimensional, the boundary condition
$\Gamma^-(0)$ is fixed at a half zero-sphere (a single point) within $E^{u}(\tilde{p}^-)$,
a distance $\varepsilon^-$ from $\tilde{p}^-$.
Similarly, if $E^{s}(\tilde{p}^+)$  is one-dimensional, the boundary condition
$\Gamma^+(0)$ is fixed at a single point within $E^{s}(\tilde{p}^+)$, a distance $\varepsilon^+$ from $\tilde{p}^+$.
However, in the case where the eigenspace $E^u(\tilde{p}^-)$ is two-dimensional,
the boundary condition
$\Gamma^-(0)$ is contained in a half one-sphere  (a half circle) of radius $\varepsilon^-$ about $\tilde{p}^-$  within
$E^u(\tilde{p}^-)$. 
Similarly,  if $E^{s}(\tilde{p}^+)$  is two-dimensional, the boundary condition
$\Gamma^+(0)$ is contained in a half circle of radius $\varepsilon^+$ about $\tilde{p}^+$  within
$E^s(\tilde{p}^+)$. 
Thus, when $E^u(\tilde{p}^-)$ or $E^s(\tilde{p}^+)$  is two-dimensional, we parameterise the boundary condition 
$\Gamma^-(0)$ on the respective half circle  by an angle parameter $\theta^-\in (0,\pi)$, and
$\Gamma^+(0)$ on the respective half circle by an angle parameter $\theta^+\in (0,\pi)$.

Third, we proceed to close the so-called \emph{Lin's gap}, which is defined as the Euclidean distance between the end points $\Gamma^-(Z^-)$ and $\Gamma^+(Z^+)$ in $\Sigma$:
$$
\eta = \left\| \Gamma^-(Z^-) -\Gamma^+(Z^+) \right\|.
$$
We are interested in structurally-stable (observable) pulse and front solutions, which typically correspond to codimension-zero heteroclinic orbits from $\tilde{p}^-$ to $\tilde{p}^+$, meaning that such connections persist on an open set of system parameters. 
%{\ch The blue sentence can be followed by something along the lines of: Codimension-zero heteroclinic orbits in a three-dimensional phase space occur when dim($E^u(\tilde{p}^-)$)+dim($E^s(\tilde{p}^+)$)=4.{\sw(Cris, I would skip the red part as I cannot see why the previous sentence is true. Can you prove it or reference the relevant theorem? I may have a counterexample.)}
%We are particularly interested in the case dim($E^u(\tilde{p}^-)$)=dim($E^s(\tilde{p}^+)$)=2.[We could discuss this when we meet]}
Hence, we close the Lin's gap by varying 
parameterized boundary conditions rather than system parameters. For example, when $W^u(\tilde{p}^-)$
and $W^s(\tilde{p}^+)$ are both two dimensional, their transverse intersections are codimension-zero heteroclinic orbits, and we
close the Lin's gap by varying
$\theta^-$ and $\theta^+$.
The strategy is to fix $\varepsilon^-$ and $\varepsilon^+$,
%in accordance with~\eqref{eq:heteroBC}, 
solve the BVP~\eqref{eq:compactifiedODE} and~\eqref{eq:leftBC} with a suitable choice of $\theta^-$, solve the  BVP~\eqref{eq:compactifiedODE} and~\eqref{eq:rightBC} with a suitable choice of $\theta^+$, and then vary $\theta^-$ and $\theta^+$ simultaneously and monitor the Lin's gap $\eta$.
Once the Lin's gap is closed, meaning that $\eta =0$ and $\Gamma^-(Z^-)=~\Gamma^+(Z^+)$,
we concatenate $\Gamma^-$ and $\Gamma^+$ to obtain a single orbit segment $\Gamma$ that satisfies the desired BVP~\eqref{eq:compactifiedODE} and \eqref{eq:heteroBC} and has a finite integration time of $Z=Z^- - Z^+$.
In this way, we can approximate a single heteroclinic orbit from
$\tilde{p}^-$ to $\tilde{p}^+$ in~\eqref{eq:compactifiedODE}, which corresponds to a pulse or front solution in~\eqref{eq:main} or~\eqref{eq:mainmf1}.

In the case where both eigenspaces $E^u(\tilde{p}^-)$ and $E^s(\tilde{p}^+)$ are of dimension two, there can be multiple coexisting heteroclinic orbits from
$\tilde{p}^-$ to $\tilde{p}^+$.
To capture multiple coexisting heteroclinic orbits, we proceed as follows.
We write $\Gamma^-_{\theta^-}$ and $\Gamma^+_{\theta^+}$ to indicate the dependence of the orbit segments $\Gamma^-$ and $\Gamma^+$ on the angle parameters $\theta^-$ and $\theta^+$.
Numerical continuation in $\theta^-$ of an orbit segment $\Gamma^-_{\theta^-}$ satisfying~\eqref{eq:leftBC} gives a parameterized family of orbit segments
$$
W^u(\tilde{p}^-) := 
\{\Gamma^-_{\theta^-}~:~\theta^-\in(0, \pi)\},
$$
which 
approximates the two-dimensional local unstable manifold of $\tilde{p}^-$.
Similarly, numerical continuation in $\theta^+$ of an orbit segment $\Gamma^+_{\theta^+}$ satisfying~\eqref{eq:rightBC} gives a parameterized family of orbit segments
$$
W^s(\tilde{p}^+) := 
\{\Gamma^+_{\theta^+}~:~\theta^+\in(0, \pi)\},
$$
which approximates the two-dimensional local stable manifold of $\tilde{p}^+$. 
$W^u(\tilde{p}^-)$ and $W^s(\tilde{p}^+)$ each intersects the cross section $\Sigma$ along a different curve.
Typically, these two curves intersect each other at isolated points in $\Sigma$. Each such isolated point in  $\Sigma$ approximates an intersection of $\Sigma$ with a different heteroclinic orbit from $\tilde{p}^-$ to $\tilde{p}^+$.
For each such isolated point in $\Sigma$, a  finite-time approximation $\Gamma$ to the corresponding heteroclinic orbit is obtained by  a concatenation of the two orbit segments, $\Gamma^-$ and $\Gamma^+$, that meet at this point.

Furthermore, in the case of multiple coexisting heteroclinic orbits from $\tilde{p}^-$ to $\tilde{p}^+$, there is a possibility that some connections become degenerate, for example, along a (codimension-one) tangency of $W^u(\tilde{p}^-)$ and $W^s(\tilde{p}^+)$, as the system parameters are varied. Such degeneracies of heteroclinic orbits in the compactified system~\eqref{eq:compactifiedODE} correspond to bifurcations of pulse and front solutions in a reaction-diffusion system~\eqref{eq:main} or~\eqref{eq:mainmf1}. 
Bifurcations of travelling waves is an area of great interest; see, for example, \cite{hagberg1994pattern,sandstede2002stability,sandstede1999essential}.
A particular advantage of our framework is that
these bifurcations can be detected by numerical continuation of the orbit segment $\Gamma$ in one of the system parameters in AUTO.  Once detected, these bifurcations can be then continued in two system parameters to produce two-parameter bifurcation diagrams of pulse and front solutions.

\section{ Pulse solutions and their bifurcations in the habitat model}
\label{sec:mhresults}
In this section, we study the existence and bifurcations of  pulse solutions in the  geographically shifting habitat model~\eqref{eq:main}--\eqref{eq:DirichletBCs}, or equivalently~\eqref{eq:mainmf1}--\eqref{eq:mainmf2}, with 
\begin{equation}
\label{eq:bchabitat}
u^\pm = U^\pm =0,   
\end{equation}
the Allee reaction term~\eqref{eq:AlleeFunction}, and the habitat function~\eqref{eq:habitat}.
%In this section, we illustrate our general approach from~\cref{sec:GenApproach} using the moving habitat model~\eqref{eq:main}--\eqref{eq:DirichletBCs} with the reaction term~\eqref{eq:AlleeFunction} and the habitat function~\eqref{eq:habitat}.
%
Specifically, we use the framework outlined in~\cref{sec:GenApproach} to obtain pulse solutions for the habitat model by computing  heteroclinic orbits in the nondimensionalized compactified system
\begin{equation}
\label{eq:compactifiedODERepeated}
  \left\{
    \begin{array}{rrl}
           \mathcal{U}_z & = & \mathcal{V}, \\
           \mathcal{V}_z & = & - \tilde{c} \, \mathcal{V} - \tilde{f}_A \Big(\mathcal{U},\tilde{H}_\alpha(s)\Big), \\
           s_z & = & \dfrac{\alpha}{2} (1-s^2).
    \end{array} \right.
\end{equation} 
Here, the dimensionless Allee reaction term~\eqref{eq:fAr} is given in terms of $s$  instead of $z$,
\begin{equation}
\label{eq:fArRepeated}
\tilde{f}_A\Big(\mathcal{U},\tilde{H}_\alpha(s)\Big) =
-\tilde{\beta}^2 \mathcal{U}
+ 4 \tilde{\beta} \, \tilde{H}_\alpha(s)\, \mathcal{U}^2
- \mathcal{U}^3,
\end{equation}
using the rescaled and extended habitat function $\tilde{H}_\alpha(s)$ from~\eqref{eq:Hcompact} with $ h^\pm =0$ and
\begin{align}
\label{eq:habitatComp}
\tilde{H}(g^{-1}_{\alpha}(s)) = \frac{e^{\tilde{L}}-e^{-\tilde{L}}}{\tanh\left(\frac{\tilde{L}}{2}\right)\left(e^{\tilde{L}}\ +e^{-\tilde{L}} + \left(\frac{1+s}{1-s}\right)^{\frac{2}{\alpha}}+\left(\frac{1-s}{1+s}\right)^{\frac{2}{\alpha}}\right)},
\end{align}
which is derived by employing the inverse coordinate transformation~\eqref{eq:inverseCF} to express $z$ in terms of $s$ in the rescaled habitat function \eqref{eq:habitatresc}.
The compactified logistic model is given in \cref{append:compactLogistic}.

\subsection{Equilibria and their stability in the  compactified habitat model
}
\label{sec:stabilityOfEq}
 We consider the two extinction equilibria for the compactified system~\eqref{eq:compactifiedODERepeated}, namely,
$$
\tilde{p}^-=(\mathcal{U},\mathcal{V},s)=(0,0,-1)\in S^-,
$$
and
$$
\tilde{p}^+=(\mathcal{U},\mathcal{V},s)=(0,0,1)\in S^+,
$$
which correspond to $u^\pm=0$ in~\eqref{eq:DirichletBCs}, or equivalently to $U^\pm =0$ in~\eqref{eq:mainmf2}.
Next, we note that  the rescaled habitat function $\tilde{H}(z)$ in~\eqref{eq:habitatresc}  decays exponentially to zero
with a decay coefficient $\rho=2$ as $z\to \pm \infty$.\footnote{Equivalently,  $\tilde{H}(g^{-1}_{\alpha}(s))$
decays to $h^\pm=0$ as $s\to\pm 1$.}
Hence, we need to choose the compactification parameter 
$$
\alpha\in(0,2],
$$
to ensure that the compactified system~\eqref{eq:compactifiedODERepeated} is C$^1$-smooth (continuously differentiable) at the added invariant planes $S^-$
containing $\tilde{p}^-$ and $S^+$ containing $\tilde{p}^+$; see section~\ref{sec:compact} and the references therein.
Linear stability analysis shows that $\tilde{p}^-$ is a saddle with eigenvalues

$$
l^-_1 = \frac{-\tilde{c} - \sqrt{\tilde{c}^2 + 4 \tilde{\beta}^2}}{2} < 0,\quad l^-_2 = \frac{-\tilde{c}+ \sqrt{\tilde{c}^2 + 4 \tilde{\beta}^2}}{2} >0,\quad\mbox{and}\quad l^-_3 = \alpha>0,
$$
 meaning that it has a two-dimensional local unstable invariant manifold $W^u(\tilde{p}^-)$. Similarly, 
$\tilde{p}^+$ is a saddle with eigenvalues 
$$
l^+_1 = l^-_1 < 0,\quad 
l^+_2 = l^-_2 > 0,
\quad\mbox{and}\quad l^+_3 = -l^-_3=-\alpha<0,
$$
meaning that it has a two-dimensional local stable invariant manifold $W^s(\tilde{p}^+)$.
Next, we note that $l_2^- \le |l_1^+|$ and, for practical convenience, we limit the choices for the compactification parameter to
$$
\alpha\in\left(0,\min\{2,l^-_2,|l_1^+|\}\right),\quad\text{or equivalently}\quad \alpha\in \left(0,\min\{2,l^-_2\}\right),
$$
to ensure that 
the additional eigenvector $v^-_3$ (corresponding to $l^-_3$) is normal to  $S^-$, the additional eigenvector $v^+_3$ (corresponding to $l^+_3$) is normal to $S^+$, and typical trajectories leave  $\tilde{p}^-$  along $v^-_3$  and approach $\tilde{p}^+$ along $v^+_3$;
see section~\ref{sec:compact} and the references therein.
%Given $\tilde{\beta}=0.15$ from \cref{table:units}, we note that$l_2^-\in(0,\tilde{\beta}]$ for $\tilde{c} \ge 0$, and take
%We then note that $l_2^-\in(0,\tilde{\beta}]$ for $\tilde{c} \ge 0$, and  $\tilde{\beta} < 2$ if and only if $D > 4.5 \times 10^{-3}$.
We then note that, for $\tilde{c} \ge 0$ and  $\tilde{\beta} < 2$, $l_2^-\in(0,\tilde{\beta}]$ if and only if $D > 4.5 \times 10^{-3}$. Hence we choose
\begin{equation}
\label{eq:alpha}
    \alpha =
  \begin{cases}
    l_2^-/2, & \text{if}\;\;\; D > 4.5 \times 10^{-3},  \\
    1, & \text{if}\;\;\;  0< D < 4.5 \times 10^{-3}.
  \end{cases}
\end{equation}

 We fix all the model parameters in~\Cref{table:units}, except for $L$, $c$, and $D$, which could be varied.
Although we work with the nondimensionalized compactified system~\eqref{eq:compactifiedODERepeated}, we specify the original parameters as input parameters in all figures.

\begin{figure}[t!]
    \centering
    \includegraphics[scale=1.0]{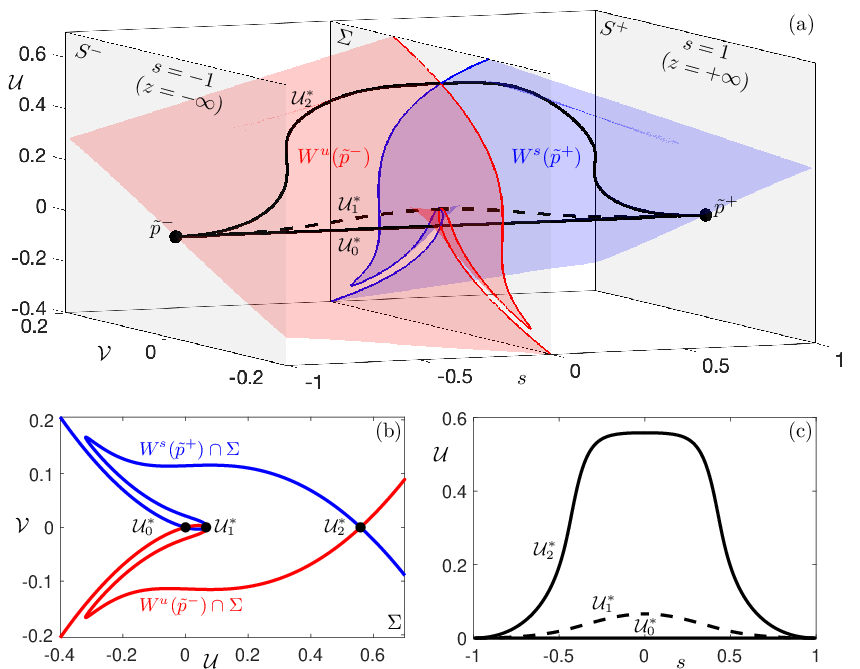}
    \caption{
    Standing-pulse solutions for the RDE~\eqref{eq:main}--\eqref{eq:DirichletBCs} with the Allee reaction term~\eqref{eq:AlleeFunction}
    obtained by computing intersections of two-dimensional invariant manifolds $W^u(\tilde{p}^-)$ and $W^s(\tilde{p}^+)$ in the compactified moving-frame ODE~\eqref{eq:compactifiedODERepeated} with the $s$-dependent Allee reaction term~\eqref{eq:fArRepeated}.
    (a) The unstable invariant manifold $W^u(\tilde{p}^-)$ of saddle $\tilde{p}^-$ (red surface), the stable invariant manifold $W^s(\tilde{p}^+)$ of saddle $\tilde{p}^+$ (blue surface), intersections of  $W^u(\tilde{p}^-)$ and $W^s(\tilde{p}^+)$  with the cross section $\Sigma$~\eqref{eq:Sigma} (red and blue curves, respectively), and transverse heteroclinic intersections of $W^u(\tilde{p}^-)$ with $W^s(\tilde{p}^+)$ (black curves $\mathcal{U}_0^*$, $\mathcal{U}_1^*$, and $\mathcal{U}_2^*$).
    (b) The cross section $\Sigma$ from (a).
    (c) The stable extinction state $\mathcal{U}_0^*$, the unstable standing pulse $\mathcal{U}_1^*$, and the stable standing pulse $\mathcal{U}_2^*$, shown over the compactified dimensionless space coordinate $s$. 
    We used $c=0 \ km/yr$, $D=~0.8\ km^2/yr$, and $L=5 \ km$; other parameter values are given in~\cref{table:units}.
    } 
    \label{fig:manifolds}
\end{figure}

\subsection{Pulse solutions as heteroclinic orbits in the compactified system} 
Codimension-zero heteroclinic orbits from $\tilde{p}^-$ to  $\tilde{p}^+$ along transverse intersections of $W^u(\tilde{p}^-)$ and $W^s(\tilde{p}^+)$ in the extended phase space $\mathbb{R}\times[-1,1]$ of
compactified system~\eqref{eq:compactifiedODERepeated} correspond to
structurally-stable pulse solutions in the habitat model, whereas codimension-one heteroclinic orbits along tangent intersections of $W^u(\tilde{p}^-)$ and $W^s(\tilde{p}^+)$ correspond to bifurcations of pulse solutions; see \cref{sec:numerical} for more details.
Numerically, we approximate both types of heteroclinic orbits using a finite-time orbit segment~\eqref{eq:orbitSegment} with boundary conditions~\eqref{eq:heteroBC} parameterized by the angles $\theta^-$ and $\theta^+$ as follows:
\begin{equation}
    \left\{
       \begin{array}{lcr}
              \Gamma(0) & = & \tilde{p}^- + \varepsilon^- (v^-_2 \cos\theta^- + v^-_3 \sin\theta^-)  , \\
               \Gamma(Z) & = & \tilde{p}^+ + \varepsilon^+ (v^+_1 \cos\theta^+ + v^+_3  \sin\theta^+).
        \end{array} \right.
       \label{eq:heteroBCexplicit}
    \end{equation}
Here, the unit eigenvectors $v^-_2$ and $v^-_3$ correspond to the eigenvalues $l^-_2$ and $l^-_3$ and span the unstable eigenspace $E^u(\tilde{p}^-)$. 
%whereas
The unit eigenvectors $v^+_1$ and $v^+_3$ correspond to the eigenvalues $l^+_1$ and $l^+_3$ and span the stable eigenspace $E^s(\tilde{p}^+)$.
To ensure that $s \in (-1,1)$ for $\Gamma(0)$ and $\Gamma(Z)$, we  use $\theta^\pm\in (0,\pi)$; see~\cref{sec:numerical}.

To  compute multiple coexisting heteroclinic orbits from $\tilde{p}^-$ to $\tilde{p}^+$, we use 
the numerical setup described in the last two paragraphs of  section~\ref{sec:numerical}.
The results shown in~\Cref{fig:manifolds} include (a)-(b) intersecting two-dimensional invariant manifolds $W^u(\tilde{p}^-)$ and 
$W^s(\tilde{p}^+)$, together with (c) the corresponding dimensionless pulse solutions $\mathcal{U}^*(s)$.
%
%The results, including intersecting two-dimensional invariant manifolds $W^u(\tilde{p}^-)$ and $W^s(\tilde{p}^+)$ together with the corresponding dimensionless pulse solutions {\ch \sout{$\mathcal{U}^*(s)$}}, are shown in~\Cref{fig:manifolds}.
%
In panel (a), we plot the (light red) unstable manifold $W^u(\tilde{p}^-)$ of $\tilde{p}^-\in S^-$ and the (light blue) stable manifold $W^s(\tilde{p}^+)$ of $\tilde{p}^+\in S^+$, each computed up to the (grey) two-dimensional cross section $\Sigma$ defined in~\eqref{eq:Sigma}. The invariant manifolds
$W^u(\tilde{p}^-)$ and $W^s(\tilde{p}^+)$ each intersect
$\Sigma$ along a different curve.
The (dark red) intersection curve $W^u(\tilde{p}^-) \cap \Sigma$ and the (dark blue) intersection curve  $W^s(\tilde{p}^+) \cap \Sigma$ are shown in more detail in~\Cref{fig:manifolds}(b).
These two curves intersect each other in three isolated (black) points in $\Sigma$, labelled $\mathcal{U}_0^*$, $\mathcal{U}_1^*$, and $\mathcal{U}_2^*$.
Each of these points corresponds to an intersection of a different heteroclinic orbit from $\tilde{p}^-$ to  $\tilde{p}^+$ with $\Sigma$.
The three heteroclinic orbits are shown (in black) in the projection onto ($s,\mathcal{U}$)-plane in panel (c).

In the habitat model~\eqref{eq:main}--\eqref{eq:DirichletBCs}, 
the trivial heteroclinic orbit $\mathcal{U}^*_0 =0$ corresponds to the extinction state, which is stable. The nontrivial heteroclinic orbits
$\mathcal{U}^*_1$ and $\mathcal{U}^*_2$ correspond to pulses  that are standing  when $c=0$, or travelling  when $c>0$; see~\Cref{fig:manifolds}(c). Pulse $\mathcal{U}^*_2$ represents the carrying capacity of the habitat and is stable. Pulse $\mathcal{U}^*_1$ is unstable and is contained in the (infinite-dimensional) {\em Allee threshold}, which separates initial states
that converge to extinction $\mathcal{U}^*_0$ from those that converge to the carrying capacity $\mathcal{U}^*_2$.\footnote{The stability of $\mathcal{U}^*_0$, $\mathcal{U}^*_1$, and $\mathcal{U}^*_2$ was obtained by numerical integration of the  habitat model~\eqref{eq:rescaledODEmoving}--\eqref{eq:habitatresc} using the method of lines~\cite{hamdi2007method,zafarullah1970application}.}

%\begin{figure}[h!]
%    \centering
%    \includegraphics[scale=1.0]{HW_habitat_f6.eps}
%    \caption{Enlargement of \cref{fig:criticalLength}(a): bifurcation diagram of the standing-wave solutions for Allee model with $c=0$ and $D=0.8$.
%    The vertical dotted lines correspond to the panels of \cref{fig:criticalLengthMechanism}.
%    } 
%    \label{fig:criticalLengthZoomed}
%\end{figure}

%\subsubsection{Critical length}
\subsection{B-tipping in a shrinking habitat: Critical length}

%{\sw\it
%1. Set $c=0$ and study one-parameter bifurcation diagrams by varying $L$.\\
%2. We obtain one-parameter bifurcation diagram of standing pulse solutions in~\eqref{} by numerical continuation of heteroclinic orbits in parameter $L$ in the compactified system~\eqref{}.
%3. The results are shown in figure4(a). Describe the figure. Identify SN bifurcation of what. Say what tipping point that is and why.
%4. Describe figure 5. Relate $L_{crit}$ to heteroclinic tangency.
%5. Compare to the logistic model. Different bif. diagram. Not B-tipping, explain why.
%}
%In \cref{sec:tippingInto}, we showed that decreasing the length of the good habitat to some critical length gives rise to B-tipping and leads to the extension of the population.{\sw (We still do not have these figures. If they are not going to materialise soon, we should leave this part out.)}

Here, we consider a static habitat with $c=0$, meaning the the original space coordinate  and the moving-frame coordinate  are identical, that is, $x=\xi$.
Thus, the new variable $s$ in~\eqref{eq:compactificationFunction} can be interpreted as
a  compactified and rescaled  original coordinate $x$:
$$
s = 
\tanh\left(\frac{1}{4}\sqrt{\frac{\beta}{D}}\, x\right).
$$

\begin{figure}[t]
    \centering
    \includegraphics[scale=1.0]{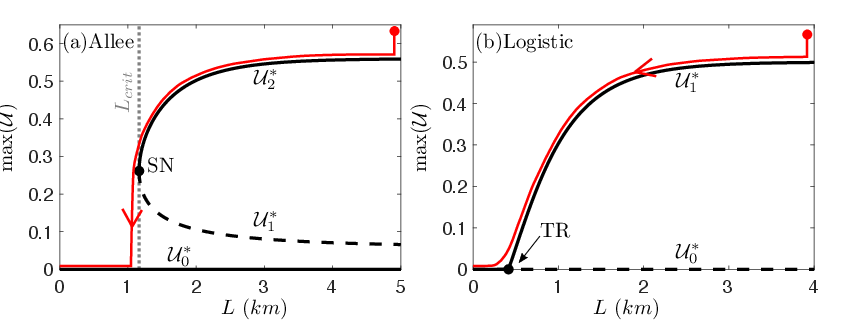}
    \caption{
    One-parameter bifurcation diagrams of standing pulses for the RDE~\eqref{eq:main}--\eqref{eq:DirichletBCs} obtained by numerical continuation of heteroclinic orbits in the compactified moving-frame ODE.
    (a) In the Allee model~\eqref{eq:compactifiedODERepeated}--\eqref{eq:fArRepeated} with $c=0 \ km/yr$ and $D=0.8 \ km^2/yr$, the branches of stable pulses $\mathcal{U}^*_2$ (solid black curve) and unstable pulses
    $\mathcal{U}^*_1$ (dashed black curve) meet and terminate in a saddle-node (SN) bifurcation (black dot), while the extinction solution $\mathcal{U}^*_0=0$ (solid black line) is stable for all values of $L$; compare with~\cref{fig:BtippingAndRtipping}(a).
    (b) In the logistic model~\eqref{eq:compactifiedODElogistic}--\eqref{eq:fLr} with $c=0 \ km/yr$ and $D=0.4 \ km^2/yr$, the branch of stable pulses $\mathcal{U}^*_1$ (solid black curve) meets the extinction solution $\mathcal{U}^*_0=0$ and terminates in a transcritical (TR) bifurcation (black dot), where $\mathcal{U}^*_0$ changes stability.
    Red trajectories are the expected solutions of the system when $L$ is decreased slowly.
    Other parameter values are given in~\cref{table:units}.
    } 
    \label{fig:criticalLength}
\end{figure}

Our aim is to describe standing-pulse solutions in the habitat model~\eqref{eq:main}--\eqref{eq:DirichletBCs} and how they depend on $L$.
To start with, we compute branches of standing pulses in~\eqref{eq:main}--\eqref{eq:DirichletBCs} by performing numerical continuation of nontrivial heteroclinic orbits in parameter $L$ in the compactified system~\eqref{eq:compactifiedODERepeated}. 
The ensuing one-parameter bifurcation diagram for the Allee reaction term~\eqref{eq:fArRepeated} is shown in~\Cref{fig:criticalLength}(a).
The stable extinction state $\mathcal{U}^*_0$ exists for all $L>0$ and is the only stationary solution for $L$ sufficiently small. 
As $L$ is increased, there is a {\em saddle-node (SN) bifurcation  of standing pulses}  at some {\em critical length $L=L_{crit}$}. 
This bifurcation gives rise to two standing pulses that exist for $L>L_{crit}$, namely, the unstable pulse $\mathcal{U}^*_1$ and the stable carrying-capacity pulse $\mathcal{U}^*_2$, and explains the the attractor diagram in Figure~\ref{fig:BtippingAndRtipping}(a).

Now consider a shrinking-habitat scenario during which  $L(t)$ slowly decreases  over time, for example, due to  deforestation and changes in land use by the growing human population.
We expect that 
the ecosystem, represented by the red trajectory in~\Cref{fig:criticalLength}(a),  tracks the stable branch of  changing carrying-capacity base states $\mathcal{U}^*_2$ until 
$L(t)$ reaches its critical value $L_{crit}$.
At this  bifurcation point the carrying-capacity base state $\mathcal{U}^*_2$ disappears. The  ensuing discontinuity in the stable branch of base states  gives rise to  a sudden transition to the alternative stable state, namely, the extinction state $\mathcal{U}^*_0$. This transition is an example of B-tipping because it is caused solely by a dangerous bifurcation of standing pulses and occurs no matter how slowly $L$ decreases.
Ecologically speaking, a habitat with $L < L_{crit}$ becomes too small to support population growth: dispersion brings the habitat population below the Allee threshold, leading to extinction.

\begin{figure}[t!]
    \centering
    \includegraphics[scale=0.9]{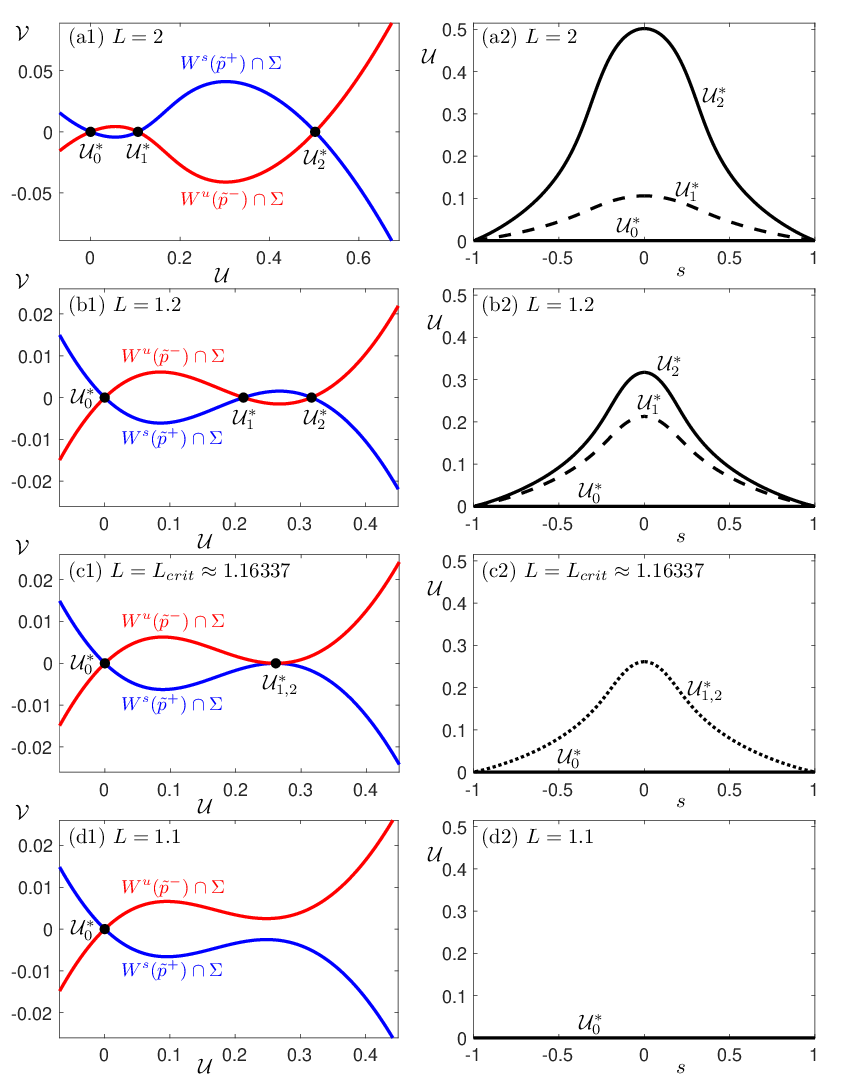}
    \caption{
   B-tipping in a shrinking habitat RDE~\eqref{eq:main}--\eqref{eq:DirichletBCs} with the Allee reaction term~\eqref{eq:AlleeFunction} shown as (left column) changing intersections of the (blue) stable  and (red) unstable  invariant manifolds on the cross section $\Sigma$~\eqref{eq:Sigma} in the compactified moving-frame ODE~\eqref{eq:compactifiedODERepeated}--\eqref{eq:fArRepeated} with $c=0 \ km/yr$ and $D=0.8 \ km^2/yr$ and (right column) colliding (solid curves) stable and (dashed curves) unstable  standing pulses, (a)--(b) before, (c) at,  and (d) after a SN bifurcation of standing pulses $\mathcal{U}^*_1$ and $\mathcal{U}^*_2$. 
   Other parameter values are given in~\cref{table:units}.
%    The geometric mechanism of B-tipping in the Allee model~\eqref{eq:compactifiedODERepeated}--\eqref{eq:fArRepeated} with $c=0 \ km/yr$ and $D=0.8 \ km^2/yr$  (a)--(b) before, (c) at, and (d) after the saddle-node bifurcation in \cref{fig:criticalLength}.
%    The left column shows the intersection curves $W^u(\tilde{p}^-) \cap \Sigma$ (red curve) and $W^s(\tilde{p}^+) \cap \Sigma$ (blue curve) and their intersections (black dots).
%    The right column shows the corresponding pulse solutions (black curves) in the ($s$,$\mathcal{U}$)-plane; solid, dashed and dotted correspond to stable, unstable and nonhyperbolic solutions, respectively.
    } 
    \label{fig:criticalLengthMechanism}
\end{figure}

The critical length $L_{crit}$ can be  detected by finding
$L$ that gives a codimension-one heteroclinic orbit along a tangent intersection of $W^u(\tilde{p}^-)$ and $W^s(\tilde{p}^+)$ in the compactified system~\eqref{eq:compactifiedODERepeated}.
This is shown in more detail in~\Cref{fig:criticalLengthMechanism}.
The left column of~\Cref{fig:criticalLengthMechanism} shows the interplay between the unstable invariant manifold $W^u(\tilde{p}^-)$ and the stable invariant manifold $W^s(\tilde{p}^+)$ on the two-dimensional cross section $\Sigma$ in the compactified system~\eqref{eq:compactifiedODERepeated}. The right column of \Cref{fig:criticalLengthMechanism} shows the rescaled stationary solutions $\mathcal{U}^*(s)$ of the habitat model~\eqref{eq:main}--\eqref{eq:DirichletBCs} that correspond to the intersections of $W^u(\tilde{p}^-)$ and $W^s(\tilde{p}^+)$.
For $L>L_{crit}$, the stable and unstable invariant manifolds intersect transversally in the three points marked with black dots; see \Cref{fig:criticalLengthMechanism}(a1) and (b1). These intersections give rise to three codimension-zero heteroclinic orbits. These orbits correspond to one trivial solution $\mathcal{U}^*_0(s)$ that exists for $L>0$,
and two standing pulses  $\mathcal{U}^*_1(s)$ and $\mathcal{U}^*_2(s)$; see \Cref{fig:criticalLengthMechanism}(a2) and (b2).
The situation is different when $L$ reaches a {\em critical level} $L=L_{crit} \approx1.16337$. In addition to the transverse intersection of the manifolds at the origin, there is a {\em tangent} intersection  away from the origin. This tangent intersection gives rise to a {\em codimension-one heteroclinic orbit}, which corresponds to a \emph{saddle-node (SN) bifurcation of standing pulses}, where 
$\mathcal{U}^*_1(s)$ and $\mathcal{U}^*_2(s)$ coalesce into $\mathcal{U}^*_{1,2}(s)$; see \Cref{fig:criticalLengthMechanism}(c1) and (c2).
For $L<L_{crit}$, there is only one transverse intersection of the manifolds at the origin, meaning that $\mathcal{U}^*_{0}(s)$ is the only stationary solution for the habitat model; see \Cref{fig:criticalLengthMechanism}(d1) and (d2).

For comparison, we show the one-parameter bifurcation diagram for the logistic growth model~\eqref{eq:compactifiedODElogistic}--\eqref{eq:fLr} in \Cref{fig:criticalLength}(b).
As $L$ is increased, there is a {\em transcritical (TR) bifurcation of standing  pulses}, in which a branch of stable carrying-capacity pulses $\mathcal{U}^*_{1}$ bifurcates from the branch of stable extinction states $\mathcal{U}^*_{0}$, while $\mathcal{U}^*_{0}$ turns unstable. The main difference from the Allee reaction term  is that TR is a safe bifurcation, meaning that there is no discontinuity in the branch of stable solutions at TR. There is no critical level $L=L_{crit}$ or bistability either. Thus, when  $L(t)$ slowly decreases over time, the (red trajectory) ecosystem tracks the branch of changing carrying-capacity base states $\mathcal{U}^*_{1}$ and declines toward the alternative extinction  state $\mathcal{U}^*_{0}$ gradually, that is, without any sudden transitions. 
In other words, there is no tipping point for the logistic reaction term.

\subsection{R-tipping in a moving habitat: Critical speed}

\begin{figure}[t]
    \centering
\includegraphics[scale=1.0]{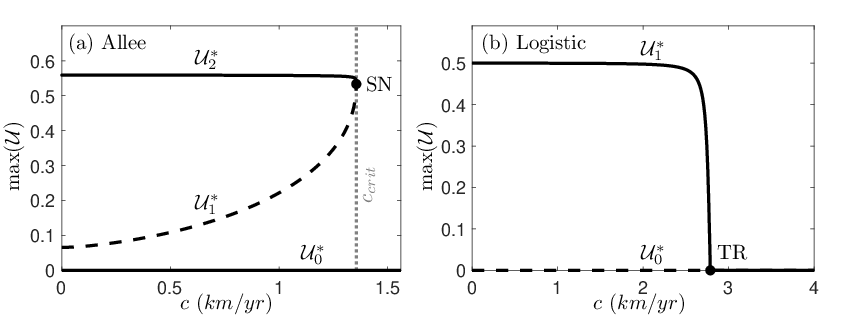}
\caption{
One-parameter bifurcation diagrams of travelling pulses for the RDE~\eqref{eq:main}--\eqref{eq:DirichletBCs} obtained by numerical continuation of heteroclinic orbits in the compactified moving-frame ODE; refer to \cref{fig:criticalLength} for more details. 
(a) The Allee model~\eqref{eq:compactifiedODERepeated}--\eqref{eq:fArRepeated} with $L=5 \ km$ and $D=0.8 \ km^2/yr$; compare with~\cref{fig:BtippingAndRtipping}(b).
(b) The logistic model~\eqref{eq:compactifiedODElogistic}--\eqref{eq:fLr} with $L=5 \ km$ and $D=0.4 \ km^2/yr$.
Other parameter values are given in~\cref{table:units}.
}   
%Bifurcation diagram of travelling-pulse solutions for the Allee model~\eqref{eq:compactifiedODERepeated}--\eqref{eq:fArRepeated} with $D=0.8 \ km^2/yr$ and $L=5 \ km$.
%Labels and colours are as for \cref{fig:criticalLength}.
\label{fig:criticalSpeed}
\end{figure}

\begin{figure}[t!]
    \centering
    \includegraphics[scale=0.9]{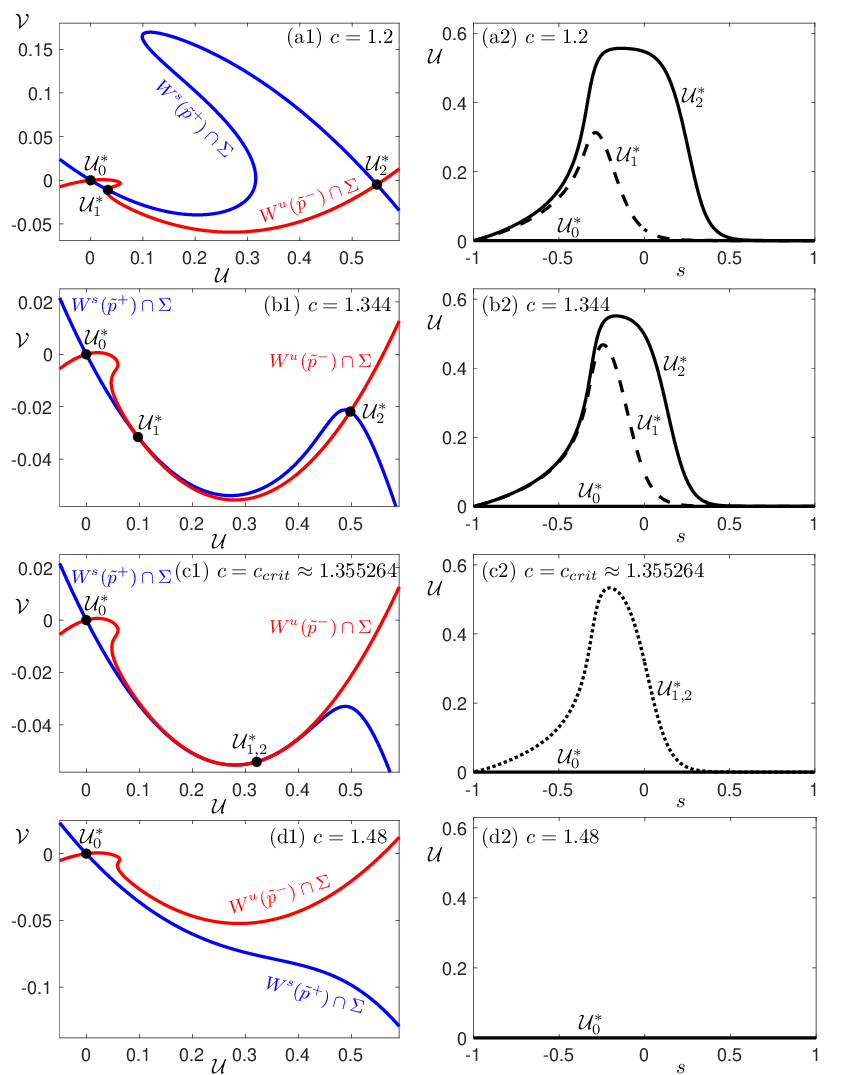}
    \caption{
     R-tipping in a shifting habitat RDE~\eqref{eq:main}--\eqref{eq:DirichletBCs} with the Allee reaction term~\eqref{eq:AlleeFunction} shown as (left column) changing intersections of the (blue) stable  and (red) unstable  invariant manifolds on the cross section $\Sigma$~\eqref{eq:Sigma} in the compactified moving-frame ODE~\eqref{eq:compactifiedODERepeated}--\eqref{eq:fArRepeated} with $D=0.8 \ km^2/yr$ and $L=5 \ km$ and (right column) colliding (solid curves) stable and (dashed curves) unstable travelling pulses, (a)--(b) before, (c) at,  and (d) after a SN bifurcation of travelling pulses. 
   Other parameter values are given in~\cref{table:units}.
%    The geometric mechanism of R-tipping in the Allee model~\eqref{eq:compactifiedODERepeated}--\eqref{eq:fArRepeated} with $D=0.8 \ km^2/yr$ and $L=5 \ km$.
%    The left column shows the intersection curves $W^u(\tilde{p}^-) \cap \Sigma$  and $W^s(\tilde{p}^+) \cap \Sigma$.
%    The right column shows the corresponding pulse solutions in the ($s$,$U$)-plane.
%     Labels and colours are as for \cref{fig:criticalLengthMechanism}.
    } 
    \label{fig:criticalSpeed2}
\end{figure}

Here, we fix $L=5 > L_{crit}$ and consider a habitat that is moving at a constant speed $c>0$, for example, due to changing weather patterns and ensuing geographical shifts in vegetation communities. Thus, the new variable $s$ in~\eqref{eq:compactificationFunction} can be interpreted as
a compactified and rescaled moving-frame coordinate $\xi=x-ct$:
$$
s = 
%\tanh\left(\frac{l_2^-}{4}\, z\right)=
\tanh \left(
\frac{-c + \sqrt{c^2 + 4\beta D}}{8D}\, (x-ct)
\right).
$$

Our aim is to describe travelling-pulse solutions in the habitat model~\eqref{eq:main}--\eqref{eq:DirichletBCs} and how they depend on $c$.
To start with, we compute branches of travelling pulses in~\eqref{eq:main}--\eqref{eq:DirichletBCs} by performing numerical continuation of nontrivial heteroclinic orbits in parameter $c$ in the compactified system. The ensuing one-parameter bifurcation diagram for the Allee model~\eqref{eq:compactifiedODERepeated}--\eqref{eq:fArRepeated} is shown in~\Cref{fig:criticalSpeed}(a).
The stable extinction state $\mathcal{U}^*_0$ exists for all $c\ge 0$.
For  $c>0$ sufficiently small, there are two travelling pulses in addition to $\mathcal{U}^*_0$, namely,  an unstable pulse $\mathcal{U}^*_1$ and a stable carrying-capacity pulse $\mathcal{U}^*_2$ that represents the ability of an ecosystem to track the moving habitat.
Interestingly, as $c$ is increased, the amplitude of the stable carrying-capacity pulse $\mathcal{U}^*_2$ remains  nearly unchanged, while the amplitude of the unstable pulse  $\mathcal{U}^*_1$ increases. Then, at some {\em critical speed} $c=c_{crit}$, there is an SN bifurcation  of travelling pulses,  at which  $\mathcal{U}^*_2$ and $\mathcal{U}^*_1$ meet and disappear. For $c > c_{crit}$, the ecosystem always goes extinct since $\mathcal{U}^*_0$ is the only stable state,  which explains the attractor diagram in Figure~\ref{fig:BtippingAndRtipping}(b).
This is an example of R-tipping because extinction is caused entirely by the rate of change in the position of the otherwise stable  ecosystem.
In other words, the spatial position of a static habitat patch in the infinite domain  has no effect on the stability of the carrying-capacity base state $\mathcal{U}^*_2$. Rather, it is the rate of change in the spatial position of the habitat patch alone that causes extinction.
Ecologically speaking, a habitat that is shifting faster than $c_{crit}$ cannot support population growth:  the dispersion  rate no longer allows the population to keep pace with the shifting  habitat, so that the population within the habitat drops below the Allee threshold, leading to extinction.

It is important to note that if the external input is a linear function of time (or, equivalently, varies at a constant speed), critical rates for R-tipping can be detected as classical autonomous bifurcations  in a suitable moving frame. This is the case here and in~\cite[Sec.3(a),(b)]{ashwin2012tipping}. However, a different approach will be required for external inputs that are nonlinear functions of time. Such inputs are left for future research.
%
%Note that, in general, critical rates for R-tipping cannot be captured by classical autonomous bifurcations. This {\ch phenomenon} may be possible for a constant speed, via a moving-frame transformation, as is the case here and in~\cite[Sec.3(a),(b)]{ashwin2012tipping}.

The critical speed $c_{crit}$ can be detected by finding
$c$ that gives a codimension-one heteroclinic orbit along a tangent intersection of $W^u(\tilde{p}^-)$ and $W^s(\tilde{p}^+)$ in the compactified system~\eqref{eq:compactifiedODERepeated}.
This is shown in more detail in~\Cref{fig:criticalSpeed2}.
The left column of~\Cref{fig:criticalSpeed2} shows the interplay between the unstable invariant manifold $W^u(\tilde{p}^-)$ and the stable invariant manifold $W^s(\tilde{p}^+)$ on the two-dimensional cross section $\Sigma$ in the compactified system~\eqref{eq:compactifiedODERepeated}. The right column of \Cref{fig:criticalSpeed2} shows the rescaled trivial and travelling-pulse  solutions $\mathcal{U}^*(s)$ of the habitat model~\eqref{eq:main}--\eqref{eq:DirichletBCs} that correspond to the intersections of $W^u(\tilde{p}^-)$ and $W^s(\tilde{p}^+)$.
For $0 < c < c_{crit}$, the stable and unstable invariant manifolds intersect transversally in the three points marked with black dots; see \Cref{fig:criticalSpeed2}(a1) and (b1). These intersections give rise to three codimension-zero heteroclinic orbits. These orbits correspond to one trivial solution $\mathcal{U}^*_0(s)$ that exists for $c>0$
and two travelling pulses  $\mathcal{U}^*_1(s)$ and $\mathcal{U}^*_2(s)$; see \Cref{fig:criticalSpeed2}(a2) and (b2). Note the increasing asymmetry in the shape of the intersecting manifolds and travelling pulses, which can be understood in terms of the symmetry-breaking advection term in~\eqref{eq:mainmf1} that is proportional to $c$.
When $c=c_{crit} \approx1.355264$, in addition to the transverse intersection of the manifolds at the origin, there is a tangent intersection  away from the origin. This tangent intersection gives rise to a codimension-one heteroclinic orbit, which corresponds to an {\em SN bifurcation of travelling pulses}, where 
$\mathcal{U}^*_1(s)$ and $\mathcal{U}^*_2(s)$ coalesce into $\mathcal{U}^*_{1,2}(s)$; see \Cref{fig:criticalSpeed2}(c1) and (c2).
For $c>c_{crit}$, there is only one transverse intersection of the manifolds at the origin, meaning that $\mathcal{U}^*_{0}(s)$ is the only stationary solution in the moving frame; see \Cref{fig:criticalSpeed2}(d1) and (d2).

For comparison, we show the one-parameter bifurcation diagram for the logistic growth model~\eqref{eq:compactifiedODElogistic}--\eqref{eq:fLr} in \Cref{fig:criticalSpeed}(b). As $c$ is increased, there is a \emph{TR bifurcation of travelling pulses} at which a branch of stable carrying-capacity pulses $\mathcal{U}^*_{1}$ 
meets the branch of unstable extinction states $\mathcal{U}^*_{0}$ and disappears, while $\mathcal{U}^*_{0}$ turns stable.
Since TR is a safe bifurcation, the stable branch of the carrying-capacity base states $\mathcal{U}^*_{1}$ declines toward the alternative extinction state $\mathcal{U}^*_{0}$ rapidly but gradually, that is, without any critical speed.
Thus, we do not consider this instability of a moving habitat with the logistic reaction term as R-tipping, but rather as a rate-induced gradual transition to extinction.

\begin{figure}[t!]
    \centering
    \includegraphics[scale=1.0]{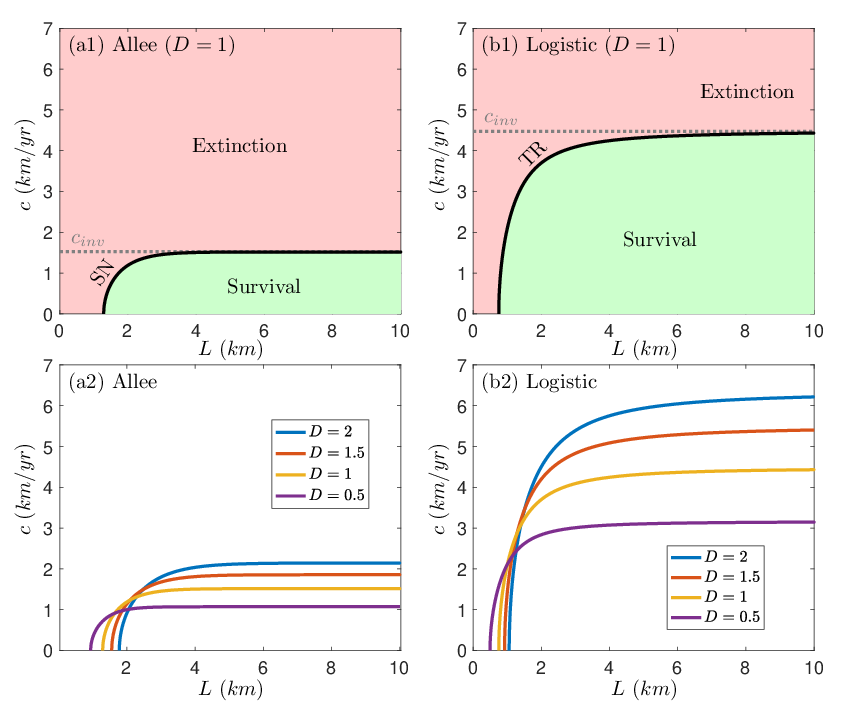} 
    \caption{
    Two-parameter bifurcation diagrams for the RDE~\eqref{eq:main}--\eqref{eq:DirichletBCs} in the parameter plane $(L,c)$ showing  the boundary between the regions of (red) extinction and (green) survival, obtained by numerical continuation of SN and TR bifurcations of heteroclinic orbits in the compactified moving-frame ODE.
    The boundary is asymptotic to (the horizontal gray line)  $c=c_{inv}$, where $c_{inv}$ depends on $D$; see \cref{appendix:invasion} for details.
    (Left column) In the Allee model~\eqref{eq:compactifiedODERepeated}--\eqref{eq:fArRepeated}, the tipping boundary is a SN bifurcation of pulses.
    (Right column) In the logistic model~\eqref{eq:compactifiedODElogistic}--\eqref{eq:fLr}, the boundary 
    is a TR bifurcation of pulses.
    (Bottom row) The boundary for different values of $D$.
    Other parameter values are given in~\cref{table:units}.
    } 
    \label{fig:bifLc}
\end{figure}

\subsection{Two-parameter bifurcation diagrams}
%Now we examine tipping regions that are triggered by both bifurcation- and rate-induced mechanisms and 
In this section, we explore two-parameter bifurcation diagrams of pulse solutions in the habitat model~\eqref{eq:main}--\eqref{eq:DirichletBCs}. %allowing three parameters $c$, $L$ and $D$ to vary.
Primarily, we are interested in the persistence of stable carrying-capacity pulses when multiple parameters are varied.
To start with, we compute the carrying-capacity pulse solution of the habitat model~\eqref{eq:main}--\eqref{eq:DirichletBCs} as a heteroclinic orbit in the compactified system and detect a codimension-one bifurcation of this pulse solution
%(e.g., saddle-node or transcritical bifurcation)
while varying a single parameter. 
Then we trace this bifurcation as a one-dimensional curve in a two-parameter plane. In this way, we identify parameter regions of {\em survival} with a stable carrying-capacity pulse and {\em extinction} with the extinction state being the only stable state.

\begin{figure}[t!]
    \centering
    \includegraphics[scale=1.0]{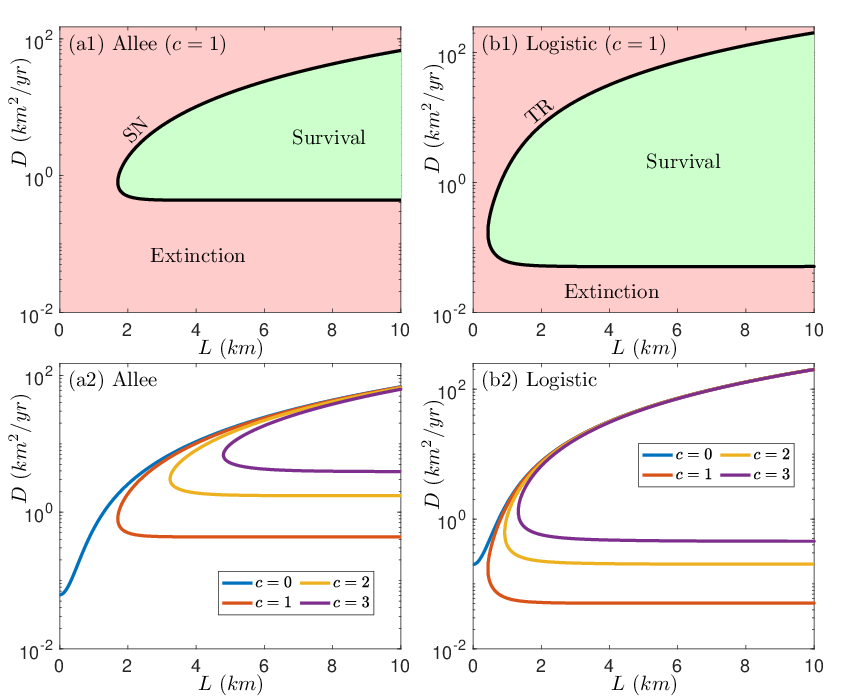}
    \caption{
    Two-parameter bifurcation diagrams for the RDE~\eqref{eq:main}--\eqref{eq:DirichletBCs} in the parameter plane $(L,D)$ showing  the boundary between the regions of (red) extinction and (green) survival, obtained by numerical continuation of SN and TR bifurcations of heteroclinic orbits in the compactified moving-frame ODE.
    See \cref{fig:bifLc} for more details.
    } 
    \label{fig:bifLD}
\end{figure}

%\Cref{fig:bifLc} shows bifurcation diagrams of pulse solutions in the ($L,c$)-plane for various values of $D$.
The two-parameter bifurcation diagram in the ($L,c$)-plane for the Allee model~\eqref{eq:compactifiedODERepeated}--\eqref{eq:fArRepeated} is shown in \Cref{fig:bifLc}(a1)--(a2).
For $D=1$, the (red) extinction  and (green) survival  regions are separated by a (black) curve of SN bifurcations of pulse solutions; see \Cref{fig:bifLc}(a1).
Note that  SN has a (dotted  grey) horizontal asymptote $c=c_{inv}$. In other words, no pulses can propagate faster than $c=c_{inv}$. It turns out that $c_{inv}$, often called the \emph{invasion speed}, is the speed of a travelling front in an infinitely-long and homogeneous habitat.\footnote{See~\cref{appendix:invasion} for more details on the computation of $c_{inv}$.}
%
%As $L \to \infty$, the separating threshold SN approaches a horizontal asymptote $c=c_{inv}$ (dotted  grey). This asymptote corresponds to the so-called \emph{invasion speed}. We define the invasion speed as the (minimum) speed of a travelling front that connects the extinction equilibrium to the carrying capacity equilibrium in the limit $L = \infty$.\footnote{For more details on the invasion speed and how we compute it see \cref{appendix:invasion}.We denote the invasion speed as $c_{inv}(D)$ to  emphasise the dependence on the dispersal rate $D$, or simply as $c_{inv}$.}
The separating tipping curves SN for different values of $D$ are shown in \Cref{fig:bifLc}(a2).
As the dispersal rate $D$ is increased, the corresponding value of the invasion speed $c_{inv}$ (not displayed) also increases, and the survival region becomes larger.
%\sout{Hence the increasing upper bound of $c$ for the survival region. Note that the survival region becomes larger when $D$ is increased.}}

For comparison, the two-parameter bifurcation diagram in the ($L,c$)-plane for the logistic growth model~\eqref{eq:compactifiedODElogistic}--\eqref{eq:fLr} is shown in \Cref{fig:bifLc}(b1)--(b2).
The extinction and survival regions for $D=1$ are separated by a TR bifurcation of pulse solutions; see \Cref{fig:bifLc}(b1).
%In the limit $L \to \infty$, the separating threshold TR approaches the horizontal asymptote $c=c_{inv}$, which corresponds to the invasion speed.
%The curve TR is representative of both critical speed and critical length.
The separating curves TR for different values of $D$ are shown in \Cref{fig:bifLc}(b2).

\begin{figure}[t!]
    \centering
    \includegraphics[scale=1.0]{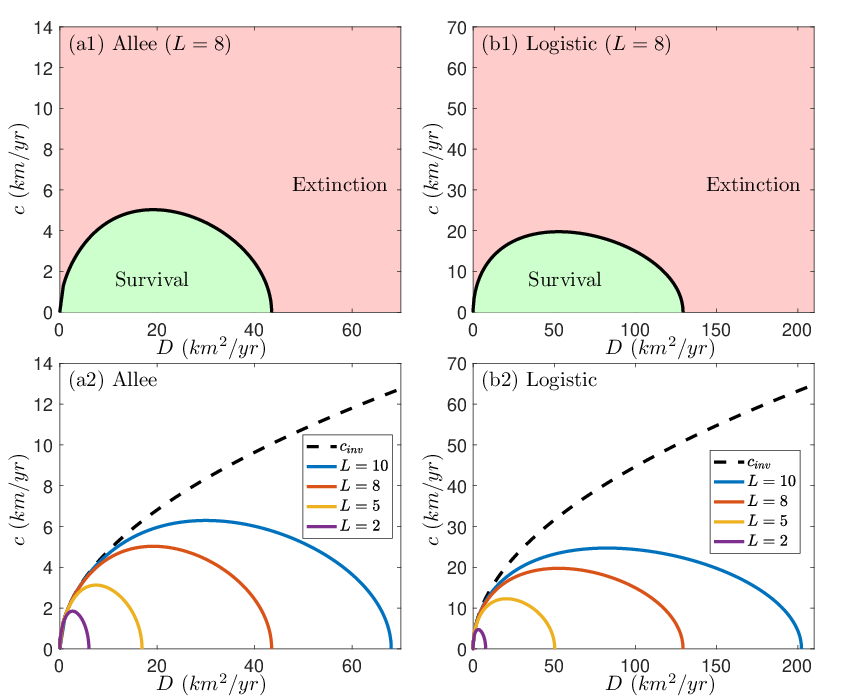}
    \caption{
    Two-parameter bifurcation diagrams for the RDE~\eqref{eq:main}--\eqref{eq:DirichletBCs} in the parameter plane $(D,c)$ showing  the boundary between the regions of (red) extinction and (green) survival, obtained by numerical continuation of SN and TR bifurcations of heteroclinic orbits in the compactified moving-frame ODE. See \cref{fig:bifLc} for more details.
    }
    \label{fig:bifDc}
\end{figure}

The two-parameter bifurcation diagram of pulse solutions in the $(L,D)$-plane for the Allee growth model~\eqref{eq:compactifiedODERepeated}--\eqref{eq:fArRepeated}  is shown in \Cref{fig:bifLD}(a1)--(a2); note the logarithmic scale of $D$.
For $c=1$, the tongue-shaped survival region is separated from extinction by a SN bifurcation of pulse solutions; see \Cref{fig:bifLD}(a1).
For sufficiently large  fixed $L$, survival is possible only within some bounded interval $(D_{min},D_{max})$ with $D_{min} > 0$.
%
%We note that dispersion is the only mechanism that allows population to move in space.
When $0< D < D_{min}$, the dispersal rate is too low for the population to keep pace with the shifting habitat, and the habitat population falls below the Allee threshold, leading to extinction.
When $D>D_{max}$,
the large dispersal rate compels a larger proportion of the population to move outside the good habitat, causes the good habitat population to drop below the Allee threshold, and also leads to extinction.
%{\ch compels the population to move outside the good habitat leading to extinction. \sout{brings the habitat population below the Allee threshold leading to extinction.}}
The observation that the population of moving habitats can only survive within a finite range of a dispersal rate has also been reported in \cite{roques2008population}.
The separating tipping curves SN for different values of $D$ are shown in \Cref{fig:bifLD}(a2). 
For $c>0$, the survival region of stable carrying-capacity travelling pulses sustain the tongue-like shape with $D_{min}>0$.
However, for $c=0$, the survival region changes shape qualitatively so that it extends to $L=0$ and retains $D_{max} >0$, but $D_{min}$ becomes zero. 
%\sout{of standing pulses is extended to $c=0$ and no longer takes a parabolic shape. Hence, for $c=0$ and fixed $L$, the survival region of standing pulses acquires some upper limit $D_{max}$ but the lower limit $D_{min}$ is always zero.}}

For comparison, the two-parameter bifurcation diagram of pulse solutions in the $(L,D)$-plane for the logistic growth model~\eqref{eq:compactifiedODElogistic}--\eqref{eq:fLr} is shown in \Cref{fig:bifLD}(b1)--(b2). %displays qualitative structure as that of the Allee growth model.

The two-parameter bifurcation diagram of pulse solutions in the $(D,c)$-plane for the Allee growth model~\eqref{eq:compactifiedODERepeated}--\eqref{eq:fArRepeated} is shown in \Cref{fig:bifDc}(a1)--(a2).
For $L=8$, the bubble-shaped survival region is separated from extinction by a SN bifurcation of pulse solutions; see \Cref{fig:bifDc}(a1).
For sufficiently small fixed $c$, the survival region exists within a bounded interval $(D_{min},D_{max})$ with $D_{min}>0$ for reasons similar to those explained in the two paragraphs above.
The separating tipping curves SN for different values of $L$ are shown in different colors in \Cref{fig:bifDc}(a2).
When $L$ is increased, the survival region extends over a wider range in the $(D,c)$-plane. 
The different tipping curves accumulate on the dashed black curve, which shows the invasion speed $c_{inv}$ for an infinitely-long homogeneous good habitat $(L=\infty)$ as a function of the dispersal rate $D$.

For comparison, the two-parameter bifurcation diagram of pulse solutions in the $(D,c)$-plane for the logistic growth model~\eqref{eq:compactifiedODElogistic}--\eqref{eq:fLr} is shown in \Cref{fig:bifDc}(b1)--(b2).

\section{Discussion}
\label{sec:discussion}

 {\bf A brief overview of the article.}
Tipping points, or critical transitions, have been studied predominantly in ordinary differential equation (ODE) models. However, they remain largely unexplored in partial differential equation (PDE) models, where spatial dynamics can give rise to new tipping mechanisms.
In this article, we studied tipping points in a special class of PDEs, namely, reaction-diffusion equations (RDEs) with a linearly time-varying and asymptotically-homogeneous reaction  term. To analyze such problems, we introduced a mathematical framework that
is based on two primary ingredients: (i) compactification of the moving-frame coordinate and (ii) computations of pulse and front solutions for such an RDE  as heteroclinic orbits connecting two equilibria from infinity in the ensuing compactified ODE.
As an illustrative example, we considered a conceptual ecosystem model subject to a geographically moving or shrinking habitat induced by climate change or human activity.
Our focus was on tipping points to extinction for the population growth model with an Allee effect (cubic nonlinearity) and how it contrasts to the simpler logistic growth model (quadratic nonlinearity).

{\bf Summary of the framework.}
The summary of our framework is as follows.
We started with a nondimensionalization of variables and parameters, followed by a reformulation of the nondimensionalized RDE into a first-order moving-frame ODE. Of importance is the fact that this ODE is only asymptotically autonomous.
Thus, we applied a compactification technique, adapted from \cite{wieczorek2021compactification}, to transform the nonautonomous ODE into an autonomous ODE on a suitably extended and compactified phase space that contains equilibria of autonomous limit systems from infinity.
In the last step of the framework, we implemented a numerical method for obtaining pulse and front solutions for the RDE by computing heteroclinic orbits connecting these equilibria in the autonomous compactified ODE.
Such heteroclinic orbits can be detected as intersections of the corresponding stable and unstable invariant manifolds using Lin's method~\cite{lin1990using}.
A particular advantage of our framework is that it also allows for numerical continuation of pulse and front solutions, as well as their bifurcations, in  the space of the system and input parameters.

{\bf Summary of the example.}
To demonstrate its applicability, the mathematical framework was implemented in AUTO~\cite{Doedel2007AUTO} and illustrated by the example of an ecosystem subject to a moving or shrinking habitat. As a result, we provided new insight into nonlinear dynamics of the moving habitat problem by performing the following:
\begin{itemize}
    \item
    Computations of the two-dimensional unstable manifold of a hyperbolic saddle from negative infinity and the two-dimensional stable invariant manifold of a hyperbolic saddle from positive infinity in the compactified autonomous ODE.
    \item
    Computation of multiple heteroclinic orbits connecting these two saddles along intersections of the manifolds. These heteroclinic orbits correspond to coexisting pulse solutions for the ecosystem RDE, and may not be possible to obtain using traditional computational techniques.
    \item
    Continuation of these heteroclinic orbits to detect their bifurcations, which correspond to bifurcations of pulse solutions for the ecosystem RDE. We distinguished between dangerous bifurcations that give rise to tipping points (abrupt transitions to extinction) and safe bifurcations that give rise to gradual transitions to extinction
    \item
    Continuation of bifurcations of heteroclinic orbits for the compactified autonomous ODE to obtain two-parameter bifurcation diagrams of pulse solutions for the ecosystem RDE.
    
\end{itemize}
%All computations were performed using the continuation software package AUTO~\cite{Doedel2007AUTO}.
Our main findings include bifurcation-induced tipping (B-tipping) to extinction below some critical length of a shrinking habitat and rate-induced tipping (R-tipping) to extinction above some critical speed of a moving habitat. We also showed that abrupt tipping points found for the Allee growth model were replaced by gradual transitions to extinction for the logistic growth model. Finally, we examined the impact of system and input parameters by analyzing curves that separate regions of  survival and extinction in two-parameter planes of: the habitat length and speed, the habitat length and population dispersion rate, and the population dispersion rate and habitat speed.

{\bf Future work.}
One interesting research direction for the future is to generalize the mathematical framework for tipping points in RDEs. Here, we considered reaction terms with linear time dependence, namely, $f(u,x-ct)$. Thus, we were able to simplify the original RDE to an ODE in the moving-frame coordinate $\xi=x-ct$. More generally, one will be interested in reaction terms with a nonlinear time dependence $g(rt)$, namely, $f(u,x-g(rt))$. 
However, such an RDE
%The challenge is that the original RDE 
no longer simplifies to an ODE in the moving-frame coordinate $\xi=x-g(rt)$.
Other interesting research directions include more complicated, possibly nonstationary, population dynamics within the habitat and extension to two spatial dimensions.

{
\section*{Acknowledgements} 
The authors would like to thank Christopher K.R.T. Jones, Jan Sieber, Bert Wuyts and, Hassan Alkhayuon for constructive comments. 
\section*{Funding} 
This work was supported by Laya Healthcare and the Enterprise Ireland Innovative Partnership Programme project IP20190771.
CRH acknowledges the partial funding by the Natural Environment Research Council (grant NE/W005042/1).
SW acknowledges partial support by the EvoGamesPlus Innovative Training Network funded by the European Union’s Horizon 2020 research and innovation program under the Marie Skłodowska-Curie grant agreement 955708.
}

%\clearpage
\appendix
\section{Appendix} \label{appendix}

\subsection{Compactified logistic model} \label{append:compactLogistic}
Here, we lay out the nondimensionalized compactified system for the logistic model
\begin{equation}
\label{eq:compactifiedODElogistic}
  \left\{
    \begin{array}{rrl}
           \mathcal{U}_z & = & \mathcal{V}, \\
           \mathcal{V}_z & = & - \tilde{c} \, \mathcal{V} - \tilde{f}_L \Big(\mathcal{U},\tilde{H}_\alpha(s)\Big), \\
           s_z & = & \dfrac{\alpha}{2} (1-s^2),
    \end{array} \right.
\end{equation} 
where the dimensionless logistic reaction term is expressed in terms of $s$ and given by
\begin{equation}
\label{eq:fLr}
    \tilde{f}_L \Big(\mathcal{U},\tilde{H}_\alpha(s)\Big) = \mathcal{U}
    \left(
    \tilde{\beta}^2 \left(2 \tilde{H}_\alpha(s)-1\right) - \mathcal{U}
    \right),
\end{equation}
and the $s$-dependent habitat function, $\tilde{H}_\alpha(s)$, is given by~\eqref{eq:Hcompact} and  \eqref{eq:habitatComp}.

\subsection{Invasion speed in a homogeneous good habitat}
\label{appendix:invasion}
%Here, we aim to find the value of $c$ above which travelling pulse solutions in moving habitat model~\eqref{eq:main}--\eqref{eq:DirichletBCs} cease to exist even for very large $L$.
%In other words, we aim to find the upper bound for the speed $c$ for travelling pulses in a homogeneous good habitat.
%We find that this upper bound always coincides with the onset of an \emph{invasion front} that connects the extinction equilibrium to the carrying capacity equilibrium  in the limit $L = \infty$.
%The speed of such invasion front is often called the \emph{invasion speed}.
%Taking the limit $L = \infty$ in the rescaled habitat function \eqref{eq:habitatresc} yields $\tilde{H} (z) =1$ for all $z$.
%Therefore, we consider the rescaled moving-frame ODE~\eqref{eq:1stOrderU}--\eqref{eq:1stOrderV} with $\tilde{H} (z) =1$.
%The ensuing ODE is the following autonomous system
Here, we compute the invasion speed in an infinite homogeneous habitat $(L=\infty$), which is the speed of travelling fronts for a constant $H(\xi)=1$.
Thus, the reaction terms in the
%both original-frame 
RDE~\eqref{eq:main} and the moving-frame ARDE~\eqref{eq:mainmf1} are now homogeneous and autonomous.
In the nondimensionalized setting, this corresponds to setting $\tilde{H}(z)=1$ for all $z$, and the nonautonomous moving-frame ODE \eqref{eq:1stOrderU}--\eqref{eq:1stOrderV} is now autonomous.
We obtain travelling fronts for the homogeneous RDE~\eqref{eq:main} by computing heteroclinic orbits in the autonomous moving-frame ODE~\eqref{eq:1stOrderU}--\eqref{eq:1stOrderV}. 

\subsubsection{Invasion speed in the Allee growth model}
\begin{figure}[t!]
    \centering
    \includegraphics[scale=1.0]{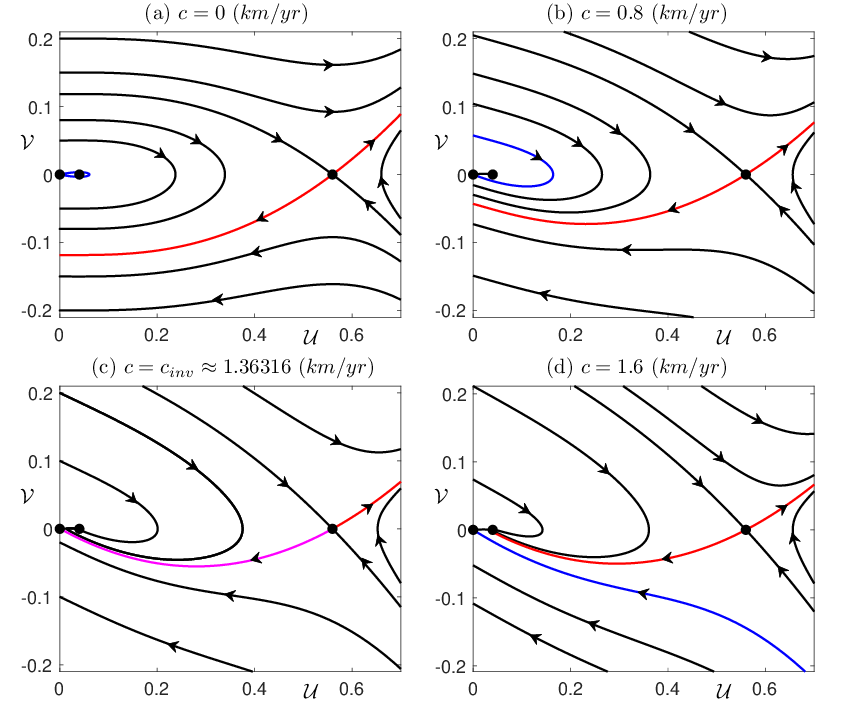}
    \caption{The phase portrait of system~\eqref{eq:invasionSpeedODE1storder}--\eqref{eq:fArInvasionS}, (a)--(b) before, (c) at, and (d) after the heteroclinic orbit. Here, $D=0.8$ and $c$ is varied.
    From left to right, the black dots on the $V$-axis are the equilibria $\mathcal{U}^*_0$ (saddle), $\mathcal{U}^*_1$ (center for $c=0$ and sink for $c>0$), and $\mathcal{U}^*_2$ (saddle).
    The blue and red orbit trajectories represent the stable and unstable manifolds of $\mathcal{U}^*_0$ and $\mathcal{U}^*_2$, respectively.
    The magenta trajectory represents the heteroclinic orbit that connects $\mathcal{U}^*_0$ with $\mathcal{U}^*_2$ at $c=c_{inv}$.
    } 
    \label{fig:invasionSpeed}
\end{figure}

For the Allee model, we consider the ensuing autonomous system
\begin{align}
\begin{split}
           \mathcal{U}_z  &=  \mathcal{V},\\
           \mathcal{V}_z & =  - \tilde{c} \, \mathcal{V} - \tilde{f}_A (\mathcal{U},1 ),
\end{split}
\label{eq:invasionSpeedODE1storder}
\end{align}
%First, we focus on the Allee growth model %\eqref{eq:fAr} to demonstrate our approach for computing the invasion speed, namely, we consider
with the Allee reaction term
\begin{equation}
\label{eq:fArInvasionS}
\tilde{f}_A(\mathcal{U},1) = \mathcal{U}
\left(
-\tilde{\beta}^2
+ 4 \tilde{\beta} \, \mathcal{U}
- \mathcal{U}^2
\right).
\end{equation}
System \eqref{eq:invasionSpeedODE1storder}--\eqref{eq:fArInvasionS} has three equilibrium points,
$$
\mathcal{U}^*_0 := (\mathcal{U},\mathcal{V}) = (0,0), \quad \mathcal{U}^*_1 := (\mathcal{U},\mathcal{V}) = ((2 - \sqrt{3})\tilde{\beta},0), \quad
\mathcal{U}^*_2 := (\mathcal{U},\mathcal{V}) = ((2 + \sqrt{3}) \tilde{\beta},0),
$$
which represent extinction, the Allee threshold, and carrying capacity, respectively.
Stability analysis shows that $\mathcal{U}^*_0$ and $\mathcal{U}^*_2$ are always of saddle type.
The Allee threshold equilibrium $\mathcal{U}^*_1$ is a center for $c=0$ and a sink for $c>0$.
To find the speed $c=c_{inv}$ of travelling fronts, we seek a heteroclinic orbit that connects the saddle equilibria $\mathcal{U}^*_0$ and $\mathcal{U}^*_2$.
We detect this heteroclinic orbit by computing the unstable manifold of $\mathcal{U}^*_2$ and the stable manifold of $\mathcal{U}^*_0$ and finding $c$ for which the two manifolds intersect; see
\Cref{fig:invasionSpeed}.
In \Cref{fig:invasionSpeed}(a)--(b), where $c<c_{inv}$ , the unstable manifold of $\mathcal{U}^*_2$ (red curve) lies below the stable manifold of $\mathcal{U}^*_0$ (blue curve).
In \Cref{fig:invasionSpeed}(c), where $c=c_{inv}$, the two manifolds intersect along a codimension-one heteroclinic orbit (magenta curve).
In \Cref{fig:invasionSpeed}(d), the unstable manifold of $\mathcal{U}^*_2$ (red curve) lies above the stable manifold of $\mathcal{U}^*_0$ (blue curve) and connects to saddle $\mathcal{U}^*_1$ instead.
Since $c_{inv}$ depends on $D$, we use numerical continuation of heteroclinic orbits
to obtain the invasion speed curve in the parameter plane $(D,c)$.
We find that the invasion speed curve $c_{inv}$ indeed marks the upper bound for the speed of travelling pulses for a bi-asymptotically homogeneous habitat; see~\cref{fig:bifLc}(a1) and~\cref{fig:bifDc}(a2).

\subsubsection{Invasion speed in the logistic growth model}
For the logistic model, the ensuing autonomous system is
\begin{align}
\begin{split}
           \mathcal{U}_z  &=  \mathcal{V},\\
           \mathcal{V}_z & =  - \tilde{c} \, \mathcal{V} - \tilde{f}_L (\mathcal{U},1 ),
\end{split}
\label{eq:invasionSpeedODElogistic}
\end{align}
with the logistic reaction term
\begin{equation}
\label{eq:fLogisticInvasionS}
\tilde{f}_A(\mathcal{U},1) =
\mathcal{U}
\left(\tilde{\beta}^2 - \mathcal{U}
\right).
\end{equation}
System \eqref{eq:invasionSpeedODElogistic}--\eqref{eq:fLogisticInvasionS} has two equilibrium points,
$$
\mathcal{U}^*_0 := (\mathcal{U},\mathcal{V}) = (0,0), \quad \mathcal{U}^*_1 := (\mathcal{U},\mathcal{V}) = (\tilde{\beta}^2,0),
$$
which represent extinction and the carrying capacity, respectively.
Stability analysis shows that $\mathcal{U}^*_0$ is always a sink and $\mathcal{U}^*_1$ is always a saddle. The invasion speed can be obtained as follows; see, for example, \cite{berestycki2009can, canosa1973nonlinear}.
When $\tilde{c}^2<4\tilde{\beta}^2$, $\mathcal{U}^*_0$ is a stable spiral,
%sink 
and when $\tilde{c}^2 \ge 4\tilde{\beta}^2$, $\mathcal{U}^*_0$ becomes a stable node.
%nodal sink. 
The heteroclinic orbit between $\mathcal{U}^*_0$ and $\mathcal{U}^*_1$ is of codimension-zero.
However, when $\tilde{c}^2 < 4\tilde{\beta}^2$, this heteroclinic orbit crosses the $\mathcal{V}$-axis  due to the spiralling nature of the spiral sink, which means that there is at least one negative $\mathcal{U}$-value along the heteroclinic orbit.  
%Hence, there is at least one negative $\mathcal{U}$-value along the heteroclinic orbit and the corresponding travelling front is physically irrelevant. 
On the other hand,  when $\tilde{c}^2 \ge 4\tilde{\beta}^2$, the heteroclinic orbit connecting $\mathcal{U}^*_0$ to $\mathcal{U}^*_1$ does not violate the condition $\mathcal{U}(z) \ge 0$ for all $z$. 
Therefore, the corresponding travelling front for $\tilde{c}^2 \ge 4\tilde{\beta}^2$ is physically relevant if and only if $\tilde{c}^2 \ge 4\tilde{\beta}^2$.
Here, we define the invasion speed in the logistic model as the value of $c$ that corresponds to the onset of physically relevant travelling fronts.
This onset is given by the analytical expression $\tilde{c}^2 = 4\tilde{\beta}^2$, or equivalently, $c=c_{inv}=2\sqrt{\beta D}$.
We find that the invasion speed $c_{inv}$ indeed marks the upper bound for the speed of travelling pulses in the logistic model; see~\cref{fig:bifLc}(b1) and~\cref{fig:bifDc}(b2).

\bibliographystyle{plain}
\bibliography{references}
\end{document}